\documentclass{article}
\usepackage{graphicx} 
\usepackage{url}
\usepackage{latexsym}
\usepackage{mismath, amsthm, latexsym, amsfonts, setspace, multirow, indentfirst, mathrsfs, pifont, xcolor, algorithm, algorithmic,mathrsfs, mathtools, amsmath, graphics, float, rotating, caption, subcaption, comment, cancel, amssymb, booktabs, bbm, braket, enumitem}
\usepackage[title]{appendix}
\usepackage{collcell}
\makeatletter
\newcolumntype{G}{>{\collectcell\@gobble}c<{\endcollectcell}@{}}
\makeatother
\usepackage[font=footnotesize]{caption}
\usepackage[
            style=apa,
            sorting = nyt]{biblatex}
\usepackage[a4paper,
            bindingoffset=0.2in,
            left=1in,
            right=1in,
            top=1in,
            bottom=1in,
            footskip=.25in]{geometry}
\usepackage[colorlinks=true,citecolor=blue,linkcolor=blue,
            urlcolor=blue]{hyperref}

  \newtheorem{theorem}{Theorem}
  
  \newtheorem{lemma}{Lemma}

\makeatletter
\renewenvironment{proof}[1][\proofname]{\par
  \pushQED{\qed}%
  \normalfont \topsep6\p@\@plus6\p@\relax
  \trivlist
  \item[\hskip\labelsep\textbf{\textit{#1}}\@addpunct{.}]\ignorespaces
}{%
  \popQED\endtrivlist\@endpefalse
}
\makeatother

    \newtheorem{definition}{Definition}
    
    \newtheorem{assumption}{Assumption}

    \newcounter{subassumption}
    \theoremstyle{TH}
    \newtheorem{fixedassumption}[assumption]{Assumption\normalfont}

    \newenvironment{assumption*}[1]
    {\renewcommand{\theassumption}{#1}\fixedassumption}
      {\endfixedassumption}
      
      \usepackage{mathtools}  
      
          \newcommand\numberthis{\addtocounter{equation}{1}\tag{\theequation}}

          \DeclareMathOperator*{\plim}{plim}
          \makeatletter

          \usepackage{arydshln}

          \bibliography{bibby} 
          
          \title{Survey Data Integration for Distribution Function and Quantile Estimation}
          \author{
            Jeremy R.A. Flood and Sayed A. Mostafa\thanks{Email: sabdelmegeed@ncat.edu} \\
            {\small\textit{Department of Mathematics \& Statistics, North Carolina A\&T State University, Greensboro, NC, USA}}
          }
          \date{}
          
          \makeatletter
          \def\adl@drawiv#1#2#3{%
          \hskip.5\tabcolsep
          \xleaders#3{#2.5\@tempdimb #1{1}#2.5\@tempdimb}%
          #2\z@ plus1fil minus1fil\relax
          \hskip.5\tabcolsep}
        \newcommand{\cdashlinelr}[1]{%
          \noalign{\vskip\aboverulesep
            \global\let\@dashdrawstore\adl@draw
            \global\let\adl@draw\adl@drawiv}
          \cdashline{#1}
            \noalign{\global\let\adl@draw\@dashdrawstore
              \vskip\belowrulesep}}
\makeatother
          
\begin{document}
\maketitle
          
\begin{abstract}
Estimates of finite population cumulative
distribution functions (CDFs) and quantiles are critical for
policy-making, resource allocation, and public health planning. For instance, federal finance agencies may require accurate estimates of the proportion of individuals with income below the federal poverty line to determine funding eligibility, while health organizations may rely on precise quantile estimates of key health variables to guide local health interventions. Despite growing interest in survey data integration, research on the integration of probability and nonprobability samples to
estimate CDFs and quantiles remains limited. In this study, we propose a novel residual-based CDF estimator that integrates information from a probability sample with data from potentially large nonprobability samples. Our approach leverages shared covariates observed in both datasets, while the response variable is available only in the nonprobability sample. Using a semiparametric approach, we train an outcome model on the nonprobability sample and incorporate model residuals with sampling weights from the probability sample to estimate the CDF of the target variable. Based on this CDF estimator, we define a quantile estimator and introduce linearization and bootstrap methods for variance estimation of both the CDF and quantile estimators. Under certain regularity conditions, we establish the asymptotic properties, including bias and variance, of the CDF estimator. Our empirical findings support the theoretical results and demonstrate the favorable performance of the proposed estimators relative to plug-in mass imputation estimators and the naïve estimators derived from the nonprobability sample only. A real data example is presented to illustrate the proposed estimators.
\end{abstract}


\section{Introduction} \label{sec:introduction}

\noindent The estimation of finite population distribution functions and quantiles plays a central role in policy-making, resource allocation, and public health surveillance. For example, quantile-based measures, such as income percentiles, poverty thresholds, or exposure levels, guide decisions in national surveys and environmental risk assessments. Despite their practical importance, distributional estimators for finite populations remain underexplored, particularly in the context of integrating probability and nonprobability survey data, where analysts have access to two complementary data sources: a probability sample $\mathcal{A}$ with known inclusion probabilities but missing outcomes $Y$, and a nonprobability sample $\mathcal{B}$ in which both the outcome $Y$ and auxiliary variables $\boldsymbol{X}$ are observed. Let $\mathcal{U}$ denote a finite population of size $N$. The goal is to estimate the finite population cumulative distribution function (CDF)
\begin{align}\label{def:cdf}
F_N(t) = \frac{1}{N}\sum_{u \in \mathcal{U}} \mathbbm{1}(Y_u \le t),
\end{align}
and the corresponding $\alpha$th quantile
\begin{align}\label{def:quant}
T_N(\alpha) = \inf\{t: F_N(t)\ge \alpha\},
\end{align}
by efficiently combining information across the two samples $\mathcal{A}$ and $\mathcal{B}$. Here $\mathbbm{1}(\cdot)$ is the usual indicator function which returns one if the condition in its argument is satisfied and zero otherwise. 

Most existing work on this context of data integration has focused on estimating the finite population mean,
\[
\mu_N = \frac{1}{N}\sum_{u \in \mathcal{U}} Y_u,
\]
under settings where the outcome $Y$ is unavailable in $\mathcal{A}$ but auxiliary covariates $\boldsymbol{X}$ are observed in both samples. Regression-based \emph{mass imputation} methods have been widely adopted in this context. The parametric and nonparametric imputation estimators proposed by \textcite{kim2021combining} and \textcite{chen2022nonparametric}, respectively, provide a foundation for borrowing strength between the two data sources under model-assisted frameworks. These methods yield asymptotically design-consistent estimators for the population mean.

However, direct extension of such approaches to CDF and quantile estimation introduces additional challenges. The empirical CDF is a nonsmooth functional of the outcome variable, and traditional linearization arguments used for mean estimation do not directly apply. Moreover, existing methods for distributional estimation under data integration are limited. \textcite{cobo2025estimation} proposed a propensity-weighted CDF estimator for nonprobability samples, but its performance deteriorates under model misspecification—a well-documented limitation of propensity-based approaches \parencite{beaumont2021pitfalls, wang2020improving, yang2020doubly, lee2011weight}.

To address these challenges, we propose a residual-based estimator for the finite population CDF, which integrates information from $\mathcal{A}$ and $\mathcal{B}$ without relying on propensity scores. The method uses residuals from an outcome model fitted on $\mathcal{B}$ to approximate the conditional distribution of $Y$ given $\boldsymbol{X}$, and applies the resulting empirical residual CDF to evaluate the probability $\Pr(Y_i \le t \mid \boldsymbol{X}_i)$ for each unit in $\mathcal{A}$. This yields a model-assisted reconstruction of the population CDF that remains consistent for $F_N(t)$. Conceptually, our approach is related to the residual-based estimator of \textcite{chambers1986estimating}, which leverages model residuals to improve CDF estimation in a totally different context where the outcome variable is observed in a probability sample and the auxiliary variables are available for the entire finite population. Additionally, we define a corresponding quantile estimator as the empirical inverse of the residual-based CDF estimator. We derive asymptotic properties of the proposed estimators, introduce both linearization and bootstrap variance estimators, and evaluate their performance through simulation studies and a real data application.

The remainder of the paper is organized as follows. Section~\ref{sec:background} introduces notation and background on survey data integration. Section~\ref{sec:reCDF} presents the proposed estimators and their theoretical properties. Section~\ref{sec:sim} reports results from simulation studies, and Section~\ref{sec:realdata} illustrates the method using a real data application. Section~\ref{sec:conclude} concludes with a discussion of key findings and directions for future research.

\section{Background and Preliminaries}
\label{sec:background}
\noindent Let $\mathcal{U} = \{1, 2, \dots, N\}$ denote the index set for a finite population of size $N$.
Let $\mathcal{A}$ denote a probability sample of size $n_\text{A}$ drawn from $\mathcal{U}$ with observed variables $\{d, \boldsymbol{X}\}$, where $d_{i} \coloneqq \pi_{i}^{-1}$ represents the sampling weight for unit $i$ with $\pi_i=\text{Pr}(i\in \mathcal{A})$ denoting the sample inclusion probability, and $\boldsymbol{X} = \begin{bmatrix} X_1 & X_2 & \cdots & X_p \end{bmatrix}$ is a matrix of auxiliary variables.
Similarly, let $\mathcal{B}$ denote a nonprobability sample of size $n_\text{B}$ from $\mathcal{U}$ with measured variables $\{Y, \boldsymbol{X}\}$, where $Y$ is the outcome (response) variable. Note that the sampling weights $d$ are not observed in $\mathcal{B}$. This data structure is summarized in Table \ref{miss} and is widely adopted in the data integration literature (e.g., \cite{kim2021combining}; \cite{yang2020statistical}). 

\begin{table}[h!]
\centering
 \caption{\label{miss} Data structure for the probability  ($\mathcal{A}$) and nonprobability ($\mathcal{B}$) samples.}
\begin{tabular*}{\textwidth}{@{\extracolsep{\fill}}lllllll@{}}
\toprule
\multicolumn{1}{l}{\textit{Sample}}&\multicolumn{1}{c}{$d$}&\multicolumn{1}{c}{$X_1$}&\multicolumn{1}{c}{$X_2$}&\multicolumn{1}{c}
{$\cdots$}&\multicolumn{1}{c}{$X_p$}&\multicolumn{1}{c}{$Y$} \\
\midrule
$\mathcal{A}$ & \checkmark & \checkmark & \checkmark & $\cdots$ & \checkmark & $\times$ \\
$\mathcal{B}$& $\times$ & \checkmark & \checkmark & $\cdots$ & \checkmark & \checkmark \\
   \bottomrule
\end{tabular*}%
\end{table}

To incorporate the auxiliary variables in estimating $F_\text{N}(t)$, we assume that the finite population is a realization from the following superpopulation model ${\xi}$: \begin{align} Y= m(\boldsymbol{X};  \boldsymbol{\beta}) + \nu(\boldsymbol{X}) \epsilon \label{sup_mod}, \end{align} where $m(\boldsymbol{X};  \boldsymbol{\beta})=\mathbb{E}\left(Y \mid \boldsymbol{X}\right)$ is a \textit{known} regression function parameterized by an unknown  vector $\boldsymbol{\beta}$, $\nu(\cdot)$ is a known, strictly positive scale function, and $\epsilon$ an independent and identically distributed (i.i.d) error term satisfying $\mathbb{E}\left(\epsilon \mid \boldsymbol{X}\right) = \mathbb{E}\left(\epsilon\right) = 0$. If all population units were observed, $\boldsymbol{\beta}$ could be estimated by solving the finite-population estimating equation $$U(\boldsymbol{\beta}) = \frac{1}{N}\sum_{u \in \mathcal{U}}{\Big(Y_{u} - m(\boldsymbol{X}_{u}; \boldsymbol{\beta})\Big)\boldsymbol{W}\left(\boldsymbol{X}_{u}; \boldsymbol{\beta}\right)} = 0$$ for some $p$-dimensional function $\boldsymbol{W}$. The resulting estimator is denoted $\boldsymbol{\beta}_{\text{N}}$, as it requires full information from $\mathcal{U}$.
Since $Y$ is observed only in $\boldsymbol{B}$, we instead estimate $\boldsymbol{\beta}$ by solving the sample-based equation $$\widehat{U}(\boldsymbol{\beta}) = \frac{1}{n_\text{B}}\sum_{j \in \mathcal{B}}{\Big(Y_{j} - m(\boldsymbol{X}_{j}; \boldsymbol{\beta})\Big)\boldsymbol{W}\left(\boldsymbol{X}_{j}; \boldsymbol{\beta}\right)} = 0,$$ whose solution is denoted $\boldsymbol{\widehat{\beta}}$. 

Before we discuss how the estimated outcome model can be used to integrate data from the two samples for estimating the finite population CDF and quantiles, we review the definitions of three common conditions for the present data integration set-up.
          
          \begin{definition}[Positivity Condition]
          \label{posit}
          Let $\delta^{\mathcal{B}}$ denote the sample indicator for membership in $\mathcal{B}$. Positivity holds if $\Pr(\delta^{\mathcal{B}} = 1 \mid \boldsymbol{X}=\boldsymbol{x})>0$ for all $\boldsymbol{x}$ in the support of $\boldsymbol{X}$. That is, for every value $\boldsymbol{x}$, there is a positive chance for the corresponding units to be selected in sample $\mathcal{B}$.  
          \end{definition}
          
         \begin{definition}[Transportability Condition]
          \label{transport}
          Let $\delta^{\mathcal{B}}$ again denote the sample indicator for $\mathcal{B}$, and let $f(Y\mid\boldsymbol{X})$ denote the conditional distribution of $Y$ given $\boldsymbol{X}$. Transportability holds if $$f(Y\mid \boldsymbol{X}, \delta^{\mathcal{B}}= 1) = f(Y\mid \boldsymbol{X}),$$  
          that is the conditional distribution of $Y$ given $\boldsymbol{X}$ is the same in $\mathcal{B}$ as in the population.
          \end{definition}
          
          Transportability is essential for transferring the model fitted on $\mathcal{B}$ to the probability sample $\mathcal{A}$. As discussed by \textcite{kim2021combining}, a sufficient condition for transportability is the \textit{ignorability condition}.

\begin{definition}[Ignorability Condition]
\label{ignore}
Let $\delta^{\mathcal{B}}$ again denote the sample indicator for $\mathcal{B}$. The ignorability condition holds if$$\Pr(\delta^{\mathcal{B}} = 1 \mid \boldsymbol{X}, Y) = \Pr(\delta^{\mathcal{B}} = 1 \mid \boldsymbol{X}),$$
i.e., after conditioning on $\boldsymbol{X}$, selection into $\mathcal{B}$ does not depend on $Y$.
\end{definition}

Definition \ref{ignore} corresponds to the \textit{missing at random} (MAR) mechanism in the missing data literature \parencite{rubin1976inference}. Under this ignorability condition, $\boldsymbol{\widehat{\beta}}$ is a consistent estimator of $\boldsymbol{\beta}$ (e.g., \textcite{kim2021combining, little2019statistical, tsiatis2006semiparametric}). With this background and notation established, in the next section, we introduce CDF and quantile estimators that effectively integrate the data available in the two samples, $\mathcal{A}$ and $\mathcal{B}$. 

\section{Proposed Estimators and Asymptotic Results} 
\label{sec:reCDF}

\subsection{Residual-based CDF Estimator}

Recall that the finite population CDF is defined as
 $$F_\text{N}(t) = \frac{1}{N}\sum_{u \in \mathcal{U}}\mathbbm{1}\left(Y_u\le t\right).$$
If the superpopulation model in \eqref{sup_mod} holds, we can write
\begin{align*} 
\mathbb{E}\left[\sum_{u \in \mathcal{U}}\mathbbm{1}\left(Y_u\le t\right)\right]
&= \sum_{u \in \mathcal{U}}{\mathbb{E}_{\mathscr{D}}\Bigg[\mathbb{E}_{\xi}\left(\mathbbm{1}\left(\epsilon_{u} \leq \frac{t-m(\boldsymbol{X}_{u}; \boldsymbol{\beta})}{\nu(\boldsymbol{X}_{u})}\right) \ \Big| \ \boldsymbol{X}_\text{N} \right)\Bigg]} \\
&= \sum_{u \in \mathcal{U}}{\Pr\left(\epsilon \leq \frac{t - m(\boldsymbol{X}_{u}; \boldsymbol{\beta})}{\nu(\boldsymbol{X}_{u})}\right)} \\
&\coloneqq \sum_{u \in \mathcal{U}}{G\left(R_{u}(t, \boldsymbol{\beta})\right)}, \numberthis \label{E_FN}
\end{align*} where $\mathbb{E}_{\mathscr{D}}$ and $ \mathbb{E}_{\xi}$ denote the design-based and model-based expectations, respectively; $\boldsymbol{X}_\text{N}$ denotes the finite population $N \times p$ matrix of covariates; $G$ is the distribution function of the model error $\epsilon$; and $$R_{u}(t; \boldsymbol{\beta}) \coloneqq \frac{t-m(\boldsymbol{X}_{u}; \boldsymbol{\beta})}{\nu(\boldsymbol{X}_{u})}.$$ 
This representation motivates a model-assisted estimator of $F_\text{N}(t)$ that integrates information from both samples $\mathcal{A}$ and $\mathcal{B}$:
\begin{align*} \widehat{F}_\text{R}(t; \boldsymbol{\widehat{\beta}}) &= \frac{1}{N}\sum_{i \in \mathcal{A}}\pi^{-1}_{i}\widehat{G}\left(R_i(t,\boldsymbol{\widehat{\beta}})\right), \numberthis \label{reCDF1} 
\end{align*}
where 
\begin{align}
\widehat{G}\left(R_i(t,\boldsymbol{\widehat{\beta}})\right) = \frac{1}{n_\text{B}}\sum_{j \in \mathcal{B}}{\mathbbm{1}\left(\widehat{\epsilon}_{j} \leq \frac{t-m(\boldsymbol{X}_{i}; \boldsymbol{\widehat{\beta}})}{\nu(\boldsymbol{X}_{i})}\right)} \label{g_def}
\end{align} 
is the empirical CDF of the estimated residuals $\widehat{\epsilon}_{j} \coloneqq \frac{Y_j-m(\boldsymbol{X}_{j}; \boldsymbol{\widehat{\beta}})}{\nu(\boldsymbol{X}_{j})}$ obtained from fitting the model in \eqref{sup_mod} to the nonprobability sample $\mathcal{B}$. The estimator in \eqref{reCDF1} combines design-based weights from the probability sample $\mathcal{A}$ with model-based residual information learned from the nonprobability sample $\mathcal{B}$. We next study its large-sample properties.

\subsection{Asymptotic Results}
\label{asymp}
We first outline the asymptotic framework for the study and the regularity assumptions under which the asymptotic properties of $\widehat{F}_\text{R}(t; \boldsymbol{\widehat{\beta}})$ are derived. Following \textcite{isaki1982survey} and \textcite{wang2011asymptotic}, let $\mathscr{U} \coloneqq \{\mathcal{U}_N\}$ denote an increasing sequence of finite populations of size $N \to \infty$, each generated from the superpopulation model \eqref{sup_mod}. For each $\mathcal{U}_N$, let $\mathcal{A}$ and $\mathcal{B}$ denote two samples drawn according to some sampling design(s), with respective sample sizes $n_\text{A}$ and $n_\text{B}$. 


\begin{assumption}
\label{a1}
Define the total sample size $n = n_\text{A} + n_\text{B}$. The sampling fractions satisfy  
\[
 \frac{n_\text{A}}{N} \to \eta_\text{A}\in (0,1), \ \frac{n_\text{B}}{N} \to \eta_\text{B}\in (0,1), \ \frac{n}{N} \to \eta \in (0,1)
\quad \text{as} \quad n_\text{A}, n_\text{B}, N \to \infty.
\] 
\end{assumption}

\begin{assumption}\label{a4}
The following conditions hold for the inclusion probabilities and the design-based variance of certain Horvitz–Thompson (HT) mean estimators from sample $\mathcal{A}$:
\begin{enumerate}[label={\textbf{4\Alph*:}}, ref={4\Alph*}]
\item \label{a4A}%
There exists positive real constants $c_{1}$ and $c_{2}$ such that, for all $u \in \mathcal{U}$, $$c_1 \leq \frac{N\pi_{u}}{n_\text{A}} \leq c_{2}.$$ 
\item \label{a4B}%
For any $Z$ with finite $2 + \delta$ population moments (for arbitrarily small $\delta >0$), there exists a positive real constant $c_{3}$ such that
$$ \frac{1}{N^2}\sum_{u \in \mathcal{U}}{\sum_{v \in \mathcal{U}}{\left\{\frac{\pi_{u,v}}{\pi_{u}\pi_{v}}-1\right\}Z_{u}Z_{v}}} \leq \frac{c_{3}}{n_\text{A}}.
$$ 
\end{enumerate}
\end{assumption}

\begin{assumption}\label{a2}
The sampling design of sample $\mathcal{B}$, although unknown, satisfies the positivity assumption in Definition \ref{posit}.
\end{assumption}

\begin{assumption}\label{a5A}
Let $\boldsymbol{\beta}^{*} := \plim_{n_\text{B} \to \infty} \boldsymbol{\widehat{\beta}}$. We assume that $\boldsymbol{\beta}^{*}$ lies within a compact set $C_{\boldsymbol{\beta}^{*}} \in \mathbb{R}$ and $\boldsymbol{\widehat{\beta}}$ is $n^{1/2}_\text{B}$--consistent for $\boldsymbol{\beta}^{*}$, in that $\boldsymbol{\widehat{\beta}} = \boldsymbol{\beta}^{*} + O_{P}(n^{-1/2}_\text{B})$. 
\end{assumption}

\begin{assumption}\label{a3}
Let $\epsilon_{j}^{*} \coloneqq \frac{Y_{j} - m(\boldsymbol{X}_{j}; \boldsymbol{\beta}^{*})}{\nu(\boldsymbol{X}_{j})}$ for $j \in \mathcal{B}$. Then $\epsilon^{*}_{1}, \epsilon^{*}_{2}, \cdots, \epsilon^{*}_{n_\text{B}} \sim_{iid} f(\epsilon^{*})$.
\end{assumption}

\begin{assumption}\label{a5}
Define
\begin{align*}
F_\text{N}(t; \boldsymbol{\beta}^{*}) &\coloneqq \frac{1}{N}\sum_{u \in \mathcal{U}}{\frac{1}{N}\sum_{v \in \mathcal{U}}{\mathbbm{1}\left(\epsilon^{*}_{v} \leq \frac{t-m(\boldsymbol{X}_{u}; \boldsymbol{\beta}^{*})}{\nu(\boldsymbol{X}_{u})}\right)}} \\
&\coloneqq\frac{1}{N}\sum_{u \in \mathcal{U}}{G_\text{N}\left(R_{u}(t; \boldsymbol{\beta}^{*})\right)}. \numberthis \label{fn_tilde}
\end{align*} 
\begin{enumerate}[label={\textbf{6\Alph*:}}, ref={5\Alph*}]
\item \label{a5B}%
    $F_\text{N}(t; \boldsymbol{\beta}^{*}) $ converges to a smooth function, say $F(t; \boldsymbol{\beta}^{*})$, uniformly in $C_{\boldsymbol{\beta}^{*}}$. 
\item \label{a5c}%
    $F(t; \boldsymbol{\beta}^{*})$ is uniformly continuous in $\boldsymbol{\beta}^{*}$ in a neighborhood of $C_{\boldsymbol{\beta}^{*}}$ and has finite first and second derivatives.    
\end{enumerate}
\end{assumption}

\begin{assumption}
\label{a6}
For $\alpha \in (\frac{1}{4}, \frac{1}{2}]$, $s$ in some compact set $C_{s} \in \mathbb{R}$, $ \boldsymbol{\beta}^{*}_{s} \coloneqq \boldsymbol{\beta}^{*} + n^{-\alpha}_\text{B}s$, and $\epsilon^{s^{*}}_{u} \coloneqq \frac{Y_{u} - m(\boldsymbol{X}_{u}; \boldsymbol{\beta}^{*}_{s})}{\nu(\boldsymbol{X}_{u})}$, the following conditions hold:  
\begin{align*}
n^{\alpha}_\text{B}\Big|F_\text{N}(t;\boldsymbol{\beta}^{*}_{s})- F_\text{N}(t;\boldsymbol{\beta}^{*})- F(t;\boldsymbol{\beta}^{*}_{s}) + F(t; \boldsymbol{\beta}^{*}) \Big| \xrightarrow{P} 0,
\end{align*}
\begin{align*}
\mathbb{E}_{\xi}\left[\frac{1}{N}\sum_{u \in \mathcal{U}}{\frac{1}{N}\sum_{v \in \mathcal{U}}{\left|\mathbbm{1}\left(\epsilon^{s^{*}}_{v} \leq \frac{t-m(\boldsymbol{X}_{u}; \boldsymbol{\beta}^{*}_{s})}{\nu(\boldsymbol{X}_{u})}\right) - \mathbbm{1}\left(\epsilon^{*}_{v} \leq \frac{t-m(\boldsymbol{X}_{u}; \boldsymbol{\beta}^{*})}{\nu(\boldsymbol{X}_{u})}\right)\right| \ \Bigg| \ \boldsymbol{X}_\text{N}}}\right] = O\left(n^{-\alpha}_\text{B}\right).
\end{align*}
\end{assumption}


Assumption \ref{a1} establishes the asymptotic framework by allowing both $n_\text{A}$ and $n_\text{B}$ to increase with the population size $N$ while controlling their relative growth rates. Assumption \ref{a4} imposes standard regularity conditions on the inclusion probabilities for sample $\mathcal{A}$, ensuring bounded weights and asymptotic design consistency of Horvitz–Thompson estimators. Assumption \ref{a2} ensures the positivity of selection probabilities in sample $\mathcal{B}$, which is essential for identifiability of the combined model. Assumption \ref{a5A} guarantees that the estimator $\boldsymbol{\widehat{\beta}}$ converges to a pseudo–true value $\boldsymbol{\beta}^{*}$ that lies in a compact parameter space; under ignorability, $\boldsymbol{\beta}^{*}$ coincides with the true $\boldsymbol{\beta}$. Assumption \ref{a3} specifies that the residuals $\epsilon^{*}_{j}$ are i.i.d. with finite variance, providing the stochastic structure needed for model-based inference from sample $\mathcal{B}$. Assumption \ref{a5} ensures the existence and smoothness of a limiting distribution function for $F_\text{N}(t; \boldsymbol{\beta}^{*})$, while Assumption \ref{a6} controls the local variation of these finite-population quantities as the parameter $\boldsymbol{\beta}^{*}$ is perturbed. Analogous smoothness and stability conditions have been used in the literature in the asymptotic analyses of model-based and model-assisted estimators (e.g., \cite{wang2011asymptotic, chambers1986estimating}).

We now establish the asymptotic properties of the CDF estimator defined in \eqref{reCDF1}. The proofs of the following results are deferred to the Appendix. 

The following lemma ensures that the estimator $\widehat{F}_\text{R}(t; \boldsymbol{\widehat{\beta}})$ is asymptotically equivalent to the infeasible version that uses the limiting parameter $\boldsymbol{\beta}^*$.

\begin{lemma}
\label{lemma_1}
Under Assumptions \ref{a1}--\ref{a6}, 
\[
\widehat{F}_\text{R}(t; \boldsymbol{\widehat{\beta}}) -\widehat{F}_\text{R}(t; \boldsymbol{\beta}^{*})
- F(t; \boldsymbol{\widehat{\beta}}) + F(t; \boldsymbol{\beta}^{*}) = o_{P}(1)
\]
and
\[
\widehat{F}_\text{R}(t; \boldsymbol{\widehat{\beta}}) -\widehat{F}_\text{R}(t; \boldsymbol{\beta}^{*}) = o_{P}(1).
\]
\end{lemma}

Lemma~\ref{lemma_1} guarantees that the estimation error in $\boldsymbol{\widehat{\beta}}$ has a negligible effect on the asymptotic behavior of the proposed CDF estimator. Building on this result, we next characterize the bias and variance of the infeasible estimator $\widehat{F}_\text{R}(t,\boldsymbol{\beta}^{*})$.

\begin{theorem}
\label{theorem_1}
Under Assumption \ref{a3}, the bias and variance of $\widehat{F}_\text{R}(t; \boldsymbol{\beta}^{*})$ are given by
\begin{align*}
\mathrm{Bias}\!\left[\widehat{F}_\text{R}(t; \boldsymbol{\beta}^{*})\right] 
&=  \frac{1}{N}\sum_{u \in \mathcal{U}}
\Big\{G\!\left(R_{u}(t, \boldsymbol{\beta^{*}})\right) - G\!\left(R_u(t;\boldsymbol{\beta})\right)\Big\}, \\
\mathrm{Var}\!\left[\widehat{F}_\text{R}(t; \boldsymbol{\beta}^{*})\right] 
&= V_{1}+ V_{2},
\end{align*}  
where
\begin{align*}
V_{1}&= \frac{1}{N^2}\sum_{u \in \mathcal{U}} \sum_{v \in \mathcal{U}}
\left\{\frac{\pi_{u,v}}{\pi_{u}\pi_{v}} -1\right\}
G\!\left(R_{u}(t; \boldsymbol{\beta}^{*})\right) 
G\!\left(R_{v}(t; \boldsymbol{\beta}^{*})\right) \\
V_{2}&= \frac{1}{n_\text{B}N^2}\sum_{u \in \mathcal{U}}\sum_{v \in \mathcal{U}}
\frac{\pi_{u,v}}{\pi_{u}\pi_{v}}
\Big\{G\big(\min\{R_{u}(t; \boldsymbol{\beta}^{*}), R_{v}(t; \boldsymbol{\beta}^{*})\}\big) 
- G\!\left(R_{u}(t; \boldsymbol{\beta}^{*})\right) G\!\left(R_{v}(t; \boldsymbol{\beta}^{*})\right)\Big\}.
\end{align*}
\end{theorem}

Theorem~\ref{theorem_1} provides explicit expressions for the bias and variance of $\widehat{F}_\text{R}(t;\boldsymbol{\beta}^*)$, decomposing the total variance into two distinct components. The first term, $V_1$, corresponds to the sampling variance of the residual-based Horvitz--Thompson estimator constructed from the probability sample $\mathcal{A}$, whereas the second term, $V_2$, captures the additional variability induced by estimating the residual distribution from the nonprobability sample $\mathcal{B}$. As expected, $V_2$ decreases as the size of the nonprobability sample, $n_\text{B}$, increases. Notice that if the ignorability condition holds, $\boldsymbol{\beta}^*=\boldsymbol{\beta}$ and the bias of $\widehat{F}_\text{R}(t;\boldsymbol{\beta}^*)$ vanishes. The subsequent theorem extends these results to the attainable estimator that replaces the limit parameter $\boldsymbol{\beta}^*$ with its consistent estimator $\boldsymbol{\widehat{\beta}}$.

\begin{theorem}
\label{theorem_2}
Under Assumptions \ref{a1}--\ref{a6}, 
\begin{align*}
\mathbb{E}\!\left(\widehat{F}_\text{R}(t; \boldsymbol{\widehat{\beta}})\right)- \mathbb{E}\left(\widehat{F}_\text{R}(t; \boldsymbol{\beta^{*}})\right)  &\to 0  \\ 
\mathrm{Var}\!\left(\widehat{F}_\text{R}(t; \boldsymbol{\widehat{\beta}})\right)
- \mathrm{Var}\!\left(\widehat{F}_\text{R}(t; \boldsymbol{\beta^{*}})\right) &\to 0.
\end{align*}
\end{theorem}

Theorem~\ref{theorem_2} shows that both the expectation and variance of the feasible estimator converge to those of the infeasible version, implying asymptotic equivalence in bias and dispersion. 

We now turn to the estimation of the variance of $\widehat{F}_R(t;\boldsymbol{\beta}^*)$. The following theorem introduces a conceptual, though unattainable, variance estimator that serves as a theoretical benchmark.

\begin{theorem}
\label{theorem_3}
Let $\widetilde{V}(t; \boldsymbol{\beta}^{*}) = \widetilde{V}_{1} + \widetilde{V}_{2}$ denote a variance estimator (unobtainable) for $\mathrm{Var}\left(\widehat{F}_\text{R}(t; \boldsymbol{\beta}^{*})\right)$, where  \begin{align*}
\widetilde{V}_{1} &= \frac{1}{(n_\text{B} - 1)N^2}\sum_{h \in \mathcal{A}}{\sum_{i \in \mathcal{A}}{\pi^{-1}_{h,i}\left\{\frac{\pi_{h,i}}{\pi_{h}\pi_{i}}-1\right\}\left\{n_\text{B}\widehat{G}\left(R_{h}\right)\widehat{G}\left(R_{i}\right) - \widehat{G}\left(\min\left\{R_{h}, R_{i}\right\}\right)\right\}}} \\
\widetilde{V}_{2} &= \frac{1}{(n_\text{B}-1)N^2}\sum_{h \in \mathcal{A}}{\sum_{i \in \mathcal{A}}{\pi^{-1}_{h,i}\left\{\frac{\pi_{h,i}}{\pi_{h}\pi_{i}}\right\}}\left\{\widehat{G}\left(\min\left\{R_{h}, R_{i}\right\}\right) -\widehat{G}\left(R_{h}\right) \widehat{G}\left(R_{i}\right)\right\}}
\end{align*}  and $R_{k} \coloneqq R_{k}(t; \boldsymbol{\beta}^{*})$.  Under Assumption \ref{a3}, $\widetilde{V}(t; \boldsymbol{\beta}^{*})$ is unbiased for $\mathrm{Var}\left(\widehat{F}_\text{R}(t; \boldsymbol{\beta}^{*})\right)$. That is,
$$
\mathbb{E}\left(\widetilde{V}(t; \boldsymbol{\beta}^{*})\right)=  \mathrm{Var}\left(\widehat{F}_\text{R}(t; \boldsymbol{\beta}^{*})\right).
$$
\end{theorem}

The variance estimator in Theorem~\ref{theorem_3} is unobtainable since it involves $\boldsymbol{\beta}^{*}$ which is not observable in practice. We next derive an attainable version that replaces the unknown terms with their sample-based estimators.

\begin{theorem}
\label{theorem_4}
Let $\widehat{V}(t; \boldsymbol{\widehat{\beta}}) = \widehat{V}_{1} + \widehat{V}_{2}$ denote an obtainable variance estimator for $\mathrm{Var}\left(\widehat{F}_\text{R}(t; \boldsymbol{\beta}^{*})\right)$,  where
 \begin{align*}
\widehat{V}_{1} &= \frac{1}{(n_\text{B} - 1)N^2}\sum_{h \in \mathcal{A}}{\sum_{i \in \mathcal{A}}{\pi^{-1}_{h,i}\left\{\frac{\pi_{h,i}}{\pi_{h}\pi_{i}}-1\right\}\left\{n_\text{B}\widehat{G}\big(\widehat{R}_{h}\big)\widehat{G}\big(\widehat{R}_{i}\big) - \widehat{G}\left(\min\left\{\widehat{R}_{h}, \widehat{R}_{i}\right\}\right)\right\}}} \\
\widehat{V}_{2} &= \frac{1}{(n_\text{B}-1)N^2}\sum_{h \in \mathcal{A}}{\sum_{i \in \mathcal{A}}{\pi^{-1}_{h,i}\left\{\frac{\pi_{h,i}}{\pi_{h}\pi_{i}}\right\}}\left\{\widehat{G}\left(\min\left\{\widehat{R}_{h}, \widehat{R}_{i}\right\}\right) -\widehat{G}\big(\widehat{R}_{h}\big) \widehat{G}\big(\widehat{R}_{i}\big)\right\}}
\end{align*} and $\widehat{R}_{k} \coloneqq R_{k}(t; \boldsymbol{\widehat{\beta}})$. Furthermore, assume there exist $M_{1}$ and $M_{2}$ satisfying $$-\infty < M_{1} \leq \frac{\pi_{h,i} - \pi_{h}\pi_{i}}{\pi_{h,i}} \leq M_{2} < \infty$$ for all $h,i \in \mathcal{A}$. Under this condition and Assumptions \ref{a1}-\ref{a6}, 
$$\left\{\widehat{V}(t; \boldsymbol{\widehat{\beta}}) - \widetilde{V}(t; \boldsymbol{\beta}^{*})\right\} = o_{P}(1).$$
\end{theorem}

Theorem~\ref{theorem_4} demonstrates that the proposed variance estimator $\widehat{V}(t; \boldsymbol{\widehat{\beta}})$ is consistent $\mathrm{Var}\left(\widehat{F}_\text{R}(t; \boldsymbol{\beta}^{*})\right)$. By Theorem~\ref{theorem_2}, $\widehat{V}(t; \boldsymbol{\widehat{\beta}})$ is also consistent for $\mathrm{Var}\left(\widehat{F}_\text{R}(t; \boldsymbol{\widehat{\beta}})\right)$. Together, Lemma~\ref{lemma_1} and Theorems~\ref{theorem_1}--\ref{theorem_4} provide the theoretical foundations for the proposed residual-based CDF estimator which integrates data from the probability and nonprobability samples.

\subsection{Alternative CDF Estimators}

\subsubsection{Na\"{i}ve Estimator}
\label{asymp_bias_naive}
Especially when sample $\mathcal{B}$ is large, it may be tempting to estimate the population CDF using only the observed outcomes from $\mathcal{B}$. This leads to the alternative CDF estimator defined as the empirical CDF of $Y$ in sample $\mathcal{B}$:
\begin{align}
\widehat{F}_\text{B}(t) \coloneqq \frac{1}{n_\text{B}}\sum_{j \in \mathcal{B}}{\mathbbm{1}\left(Y_{j} \leq t\right)}.\label{naive_cdf}
\end{align}
Because this estimator ignores potential bias in sample $\mathcal{B}$ arising from the absence of a known sampling design (i.e., $\mathcal{B}$ is a nonprobability sample), it is referred to as the \textit{na\"{i}ve estimator}. 

We will briefly analyze the bias of $\widehat{F}_\text{B}(t)$. Recall that $\delta^{\mathcal{B}}$ denotes the sample membership indicator for $\mathcal{B}$ and note that conditional on $\delta^{\mathcal{B}}$ and $\boldsymbol{X}_\text{N}$, we have
\begin{align*}
\mathbb{E}_{\xi}\!\left(\widehat{F}_\text{B}(t) \,\big|\, \delta^{\mathcal{B}}, \boldsymbol{X}_\text{N}\right)
&= \mathbb{E}_{\xi}\!\left(
\frac{1}{n_\text{B}}\sum_{j \in \mathcal{B}} \mathbbm{1}\!\left(Y_{j} \leq t\right)
\,\Big|\, \delta^{\mathcal{B}}, \boldsymbol{X}_\text{N}\right) \\
&= \mathbb{E}_{\xi}\!\left(
\frac{1}{n_\text{B}}\sum_{j \in \mathcal{B}}
\mathbbm{1}\!\left(\epsilon_{j} \leq \frac{t - m(\boldsymbol{X}_{j}; \boldsymbol{\beta})}
{\nu(\boldsymbol{X}_{j})}\right)
\,\Big|\, \delta^{\mathcal{B}}, \boldsymbol{X}_\text{N}\right) \\
&= \frac{1}{n_\text{B}}\sum_{j \in \mathcal{B}} G\!\left(R_{j}(t; \boldsymbol{\beta})\right),
\end{align*}
under the ignorability condition for sample~$\mathcal{B}$.  
Taking expectation with respect to the sampling design~$\mathscr{D}$ gives
\begin{align*}
\mathbb{E}_{\mathscr{D}}\!\left[\mathbb{E}_{\xi}\!\left(\widehat{F}_\text{B}(t) \,\Big|\, \delta^{\mathcal{B}}, \boldsymbol{X}_\text{N}\right)\right]
&= \frac{1}{n_\text{B}} \sum_{u \in \mathcal{U}} \Pr(u \in \mathcal{B}) \, G\!\left(R_{u}(t; \boldsymbol{\beta})\right) \\
&= \frac{1}{N}\sum_{u \in \mathcal{U}} f^{-1}_{n_\text{B}}\Pr(u \in \mathcal{B}) \, G\!\left(R_{u}(t; \boldsymbol{\beta})\right),
\end{align*}
where $f_{n_\text{B}} = n_\text{B}/N$ is the sampling fraction of sample~$\mathcal{B}$.

In contrast, for the finite-population CDF,
\begin{align*}
\mathbb{E}\!\left(F_\text{N}(t)\right)
&= \mathbb{E}_{\mathscr{D}}\!\left[\mathbb{E}_{\xi}\!\left(F_\text{N}(t) \,\big|\, \boldsymbol{X}_\text{N}\right)\right] \\
&= \mathbb{E}_{\mathscr{D}}\!\left[\frac{1}{N}\sum_{u \in \mathcal{U}} G\!\left(R_{u}(t; \boldsymbol{\beta})\right)\right] \\
&= \frac{1}{N}\sum_{u \in \mathcal{U}} G\!\left(R_{u}(t; \boldsymbol{\beta})\right).
\end{align*}
Therefore, the bias of $\widehat{F}_\text{B}(t)$ relative to $F_\text{N}(t)$ is
\[
\mathbb{E}\!\left[\widehat{F}_\text{B}(t) - F_\text{N}(t)\right]
= \frac{1}{N}\sum_{u \in \mathcal{U}} G\!\left(R_{u}(t; \boldsymbol{\beta})\right)
\Big\{ f^{-1}_{n_\text{B}}\Pr(u \in \mathcal{B}) - 1 \Big\}.
\]
Hence, $\widehat{F}_\text{B}(t)$ is generally biased for $F_\text{N}(t)$ unless 
$f^{-1}_{n_\text{B}}\Pr(u \in \mathcal{B}) \approx 1$, or equivalently, the inclusion probabilities for sample $\mathcal{B}$ are proportional to the sampling fraction.  
Alternatively, if $G\!\left(R_{u}(t; \boldsymbol{\beta})\right)$ is constant across the population and 
$\sum_{u \in \mathcal{U}}\Pr(u \in \mathcal{B}) = n_\text{B}$, then
\begin{align*}
\frac{1}{N}\sum_{u\in\mathcal{U}} G\!\left(R_{u}(t;\boldsymbol{\beta})\right)
\Big\{f^{-1}_{n_\text{B}}\Pr(u\in\mathcal{B})-1\Big\}
&= c\!\left\{\frac{1}{N}\frac{N}{n_\text{B}}\sum_{u\in\mathcal{U}}\Pr(u\in\mathcal{B}) - 1\right\} \\
&= c\!\left\{\frac{1}{n_\text{B}}\,n_\text{B} - 1\right\} = 0,
\end{align*}
implying that the bias vanishes in this special case.

In practice, $G\!\left(R_{u}(t;\boldsymbol{\beta})\right)$ tends to be nearly constant for sufficiently large or small values of $t$ (corresponding to the tails of the response distribution).  
Consequently, the na\"{i}ve estimator $\widehat{F}_\text{B}(t)$ may perform reasonably well in the distributional extremes, even though it remains biased elsewhere.

\subsubsection{Plug-In Estimator}
\label{sec:plugcdf}

An alternative approach for estimating the CDF is to adopt a traditional \textit{mass imputation} strategy, where each unobserved $Y_{i}$ is replaced by its model-predicted value, $\widehat{Y}_{i} \coloneqq m(\boldsymbol{X}_{i}; \boldsymbol{\widehat{\beta}})$. Let $\widehat{F}_\text{P}(t; \boldsymbol{\widehat{\beta}})$ denote the resulting \textit{plug-in CDF estimator}, formally defined as  
\begin{align*} 
\widehat{F}_\text{P}(t; \boldsymbol{\widehat{\beta}}) 
&= \frac{1}{N}\sum_{i \in \mathcal{A}}{\pi^{-1}_{i}\,\mathbbm{1}\!\left(\widehat{Y}_{i} \leq t\right)}\\
&= \frac{1}{N}\sum_{i \in \mathcal{A}}{\pi^{-1}_{i}\,\mathbbm{1}\!\left(m(\boldsymbol{X}_{i}; \boldsymbol{\widehat{\beta}}) \leq t \right)}. \numberthis \label{plugin_cdf}
\end{align*}

Even in the hypothetical case where $\boldsymbol{\widehat{\beta}} = \boldsymbol{\beta}$, we have
\begin{align*}
\mathbb{E}_{\mathscr{D}}\Big(\widehat{F}_\text{P}(t; \boldsymbol{\widehat{\beta}})\Big|_{\boldsymbol{\widehat{\beta}} = \boldsymbol{\beta}}\Big) 
&= \frac{1}{N}\sum_{u \in \mathcal{U}}{\mathbbm{1}\!\Big(m(\boldsymbol{X}_{u}; \boldsymbol{\beta}) \leq t \Big)} \\
&= \frac{1}{N}\sum_{u \in \mathcal{U}}{\mathbbm{1}\!\Big(Y_{u} - \nu(\boldsymbol{X}_{u})\epsilon_{u} \leq t \Big)},
\end{align*}
which implies that $\widehat{F}_\text{P}(t; \boldsymbol{\widehat{\beta}})\big|_{\boldsymbol{\widehat{\beta}} = \boldsymbol{\beta}}$ is not design-unbiased for $F_\text{N}(t)$ under model~\eqref{sup_mod}, unless $\nu(\boldsymbol{X}_{u})\epsilon_{u}$ is effectively zero for all $u \in \mathcal{U}$. Furthermore, 
\begin{align*}
\mathbb{E}_{\mathscr{D}}\!\Bigg[\mathbb{E}_{\xi}\!\left(\frac{1}{N}\sum_{u \in \mathcal{U}}{\frac{\delta_{u}^{\mathcal{A}}}{\pi_{u}}\,
\mathbbm{1}\!\Big(m(\boldsymbol{X}_{u}; \boldsymbol{\beta})\leq t \Big)} \ \Big| \ \delta^{\mathcal{A}}, \boldsymbol{X}_\text{N} \right)\Bigg] 
&= \frac{1}{N}\sum_{u\in\mathcal{U}}{\Pr\!\left(m(\boldsymbol{X}_{u}; \boldsymbol{\beta}) \leq t\right)} \\
&\neq \frac{1}{N}\sum_{u \in \mathcal{U}}{G\!\left(R_{u}(t; \boldsymbol{\beta})\right)},
\end{align*}
which indicate that the bias in the plug-in estimator $\widehat{F}_\text{P}(t; \boldsymbol{\widehat{\beta}})$ remains under the combined design-model-based framework. In contrast, the proposed residual-based CDF estimator $\widehat{F}_\text{R}(t; \boldsymbol{\widehat{\beta}})$ incorporates the estimated residual distribution to correct this bias and improve the efficiency in estimating $F_\text{N}(t)$.
 
\subsection{Quantile Estimation}
\label{sec:quante}

We conclude this section with a discussion of quantile estimation. From Eq.~\eqref{def:quant}, an estimator of the $\alpha$-th finite population quantile, $T_\text{N}(\alpha)$, based on the proposed CDF estimator $\widehat{F}_\text{R}(t; \boldsymbol{\widehat{\beta}})$ is defined as 
\begin{align}
\widehat{T}_\text{R}(\alpha) = \inf\!\left\{t : \widehat{F}_\text{R}(t; \boldsymbol{\widehat{\beta}}) \geq \alpha \right\}, \label{t_r}
\end{align}
with analogous estimators obtained from $\widehat{F}_\text{P}(t; \boldsymbol{\widehat{\beta}})$ and $\widehat{F}_\text{B}(t)$ as
\begin{align}
\widehat{T}_\text{P}(\alpha) &= \inf\!\left\{t : \widehat{F}_\text{P}(t; \boldsymbol{\widehat{\beta}}) \geq \alpha \right\}, \label{t_p} \\
\widehat{T}_\text{B}(\alpha) &= \inf\!\left\{t : \widehat{F}_\text{B}(t) \geq \alpha \right\}. \label{t_naive}
\end{align}
Here, $\alpha$ denotes the percentile of interest. 

Following \textcite{kovar1988bootstrap} and \textcite{francisco1991quantile}, the variance of $\widehat{T}(\alpha)$ can be approximated by
\begin{align}
\widehat{\mathrm{Var}}\!\left[\widehat{T}(\alpha)\right] 
= \left(\frac{\widehat{T}_\text{UL} - \widehat{T}_\text{LL}}{2\,z_{\gamma/2}}\right)^{2}. \label{var_q}
\end{align}
where
\begin{align}
\widehat{T}_\text{LL} &= \inf\!\left\{t : \widehat{F}(t) \geq \alpha - z_{\gamma/2}\,\widehat{\mathrm{SE}}\!\left[\widehat{F}\!\left(\widehat{T}(\alpha)\right)\right]\right\}, \label{q:LL} \\
\widehat{T}_\text{UL} &= \inf\!\left\{t : \widehat{F}(t) \geq \alpha + z_{\gamma/2}\,\widehat{\mathrm{SE}}\!\left[\widehat{F}\!\left(\widehat{T}(\alpha)\right)\right]\right\}, \label{q:UL}
\end{align}
represent the lower and upper limits of the \textcite{woodruff1952confidence} large-sample confidence interval for quantiles under general sampling designs, with $z_{\gamma/2}$ denoting the $(1-\gamma/2)^{\text{th}}$ quantile of the standard normal distribution and $\widehat{\mathrm{SE}}$ the estimated standard error of $\widehat{F}(t)$. \textcite{sitter2001note} demonstrated that the Woodruff interval performs well even when the corresponding interval for the distribution function performs poorly. An alternative variance estimator based on the bootstrap method will be defined and empirically evaluated in the next section.

\section{Simulation Study}
\label{sec:sim}

\noindent We conducted a Monte Carlo simulation study to evaluate the finite-sample performance of the proposed CDF and quantile estimators, $\widehat{F}_\text{R}(t; \boldsymbol{\widehat{\beta}})$ and $\widehat{T}_\text{R}(\alpha)$, relative to their plug-in and na\"{i}ve counterparts, $\widehat{F}_\text{P}(t; \boldsymbol{\widehat{\beta}})$ and $\widehat{T}_\text{P}(\alpha)$, and $\widehat{F}_\text{B}(t)$ and $\widehat{T}_\text{B}(\alpha)$. 

To assess relative efficiency, we used the root mean squared error ratio (RMSER), defined for a generic estimator $\widehat{\theta}$ as 
\begin{align}
\mathrm{RMSER}\left(\widehat{\theta}\right) 
= \sqrt{\frac{\mathrm{MSE}\left(\widehat{\theta}\right)}{\mathrm{MSE}\left(\widehat{\theta}_{\pi}\right)}}, 
\label{rrmse_def}
\end{align}
where $\widehat{\theta}_{\pi}$ denotes a ``gold-standard'' design-based estimator from sample $\mathcal{A}$ if the response variable $Y$ was observed in $\mathcal{A}$. In particular, we take $\widehat{\theta}_{\pi} \coloneqq \{\widehat{F}_{\pi}(t), \widehat{T}_{\pi}(\alpha)\}$, representing the Horvitz–Thompson (HT) estimators of $F_\text{N}(t)$ and $T_\text{N}(\alpha)$, respectively.

\subsection{Simulation Settings}

To thoroughly assess the robustness and efficiency of the proposed estimators, we considered four superpopulation models, $\{\xi_{i}\}_{i=1}^{4}$, that vary in functional complexity between $Y$ and $\boldsymbol{X}$:

\begin{itemize}
\item \textbf{Model $\xi_{1}$:}
$$
Y = \sum_{a=1}^{2}{4X_{a}} + \sum_{b=3}^{4}{2X_{b}} + \epsilon,
$$
where $X_{1}, X_{2} \sim \mathrm{N}(\mu=2, \sigma=1)$, $X_{3}, X_{4} \sim \mathrm{N}(4, 1)$, and $\epsilon \sim \mathrm{N}(0, 3)$.

\item \textbf{Model $\xi_{2}$:}
$$
Y = \sum_{a=1}^{2}{4X_{a}^{2}} + \sum_{b=3}^{4}{2X_{b}^{2}} 
+ \sum_{c \in \{1,3\}}{\left(X_{c} + X_{c+1}\right)^{2}} + \epsilon,
$$
where $X_{1}, X_{2} \sim \mathrm{Uniform}(0,4)$, $X_{3}, X_{4} \sim \mathrm{Uniform}(4,8)$, and $\epsilon \sim \mathrm{N}(0, 50)$.

\item \textbf{Model $\xi_{3}$} \parencite{scornet2016random, maia2021predictive}:
$$
Y = -\sin(X_{1}) + X_{2}^{2} + X_{3} - \mathrm{e}^{-X_{4}^{2}} + \epsilon,
$$
where $X_{1}, \dots, X_{4} \sim \mathrm{Uniform}(-1, 1)$ and $\epsilon \sim \mathrm{N}(0, \sqrt{0.5})$.

\item \textbf{Model $\xi_{4}$} \parencite{roy2012robustness, maia2021predictive}:
$$
Y = X_{1} + 0.707X_{2}^{2} + 2\mathbbm{1}(X_{3}>0) 
+ 0.873\ln(|X_{1}|)|X_{3}| + 0.894X_{2}X_{4} 
+ 2\mathbbm{1}(X_{5}>0) + 0.46\mathrm{e}^{X_{6}} + \epsilon,
$$
where $X_{1}, \dots, X_{6} \sim \mathrm{N}(0, 1)$ and $\epsilon \sim \mathrm{N}(0, 1)$.
\end{itemize}

For each model, we generated a finite population of size $N = 100{,}000$. From each population, a simple random sample without replacement (SRS) of size $n_\text{A} = 1{,}000$ was drawn to represent the probability sample $\mathcal{A}$. For the nonprobability sample $\mathcal{B}$, we considered sample sizes $\boldsymbol{n}_\text{B} = [n_\text{A}, 10n_\text{A}, 20n_\text{A}]$ and generated both missing-at-random (MAR) and missing-not-at-random (MNAR) versions using stratified SRS. Specifically, to simulate MAR, we identified $X^{*}$ as the covariate with the highest Pearson correlation with $Y$. Each population element $u \in \mathcal{U}$ was assigned to stratum I if $X_{u}^{*} \leq T_\text{N}(X^{*}; \alpha = 0.5)$ and to stratum II otherwise, where $T_\text{N}(X^{*}; \alpha = 0.5)$ denotes the finite population median of $X^{*}$. We then selected $n_\text{I} = 0.15\boldsymbol{n}_\text{B}$ and $n_\text{II} = 0.85\boldsymbol{n}_\text{B}$ units from the respective strata. The MNAR case followed the same stratified procedure, except strata were defined using $Y_{u} \leq T_\text{N}(Y; \alpha = 0.5)$ instead of $X^{*}$. 

Using each of the $n_\text{sim} = 1{,}500$ randomly drawn $(\mathcal{A}, \mathcal{B})$ pairs, we estimated the CDF and quantiles at the 1st, 10th, 25th, 50th, 75th, 90th, and 99th percentiles using the following estimators:

\begin{enumerate}
\item \textbf{Na\"{i}ve Estimators:}
\begin{itemize}
\item $\widehat{F}_\text{B}(t)$ defined in Eq.~\eqref{naive_cdf}
\item $\widehat{T}_\text{B}(\alpha)$ defined in Eq.~\eqref{t_naive}
\end{itemize}

\item \textbf{Plug-in Estimators:}
\begin{itemize}
\item $\widehat{F}_\text{P}(t; \boldsymbol{\widehat{\beta}})$ defined in Eq.~\eqref{plugin_cdf}
\item $\widehat{T}_\text{P}(\alpha)$ defined in Eq.~\eqref{t_p}
\end{itemize}

\item \textbf{Proposed Residual-based Estimators:}
\begin{itemize}
\item $\widehat{F}_\text{R}(t; \boldsymbol{\widehat{\beta}})$ defined in Eq.~\eqref{reCDF1}
\item $\widehat{T}_\text{R}(\alpha)$ defined in Eq.~\eqref{t_r}
\end{itemize}
\end{enumerate}

Both $\widehat{F}_\text{P}(t; \boldsymbol{\widehat{\beta}})$ and $\widehat{F}_\text{R}(t; \boldsymbol{\widehat{\beta}})$ employed multivariable linear regression, fitted on sample $\mathcal{B}$ using the \texttt{lm()} function in \textsf{R} \parencite{RcoreTeam}. All computations were performed in \textsf{R}, and the full simulation code is available in the arXiv preprint by \textcite{flood2024survey}.

\begin{figure}[h!]
  \centering
  \includegraphics[width=0.9\textwidth,height=0.9\textheight,keepaspectratio]{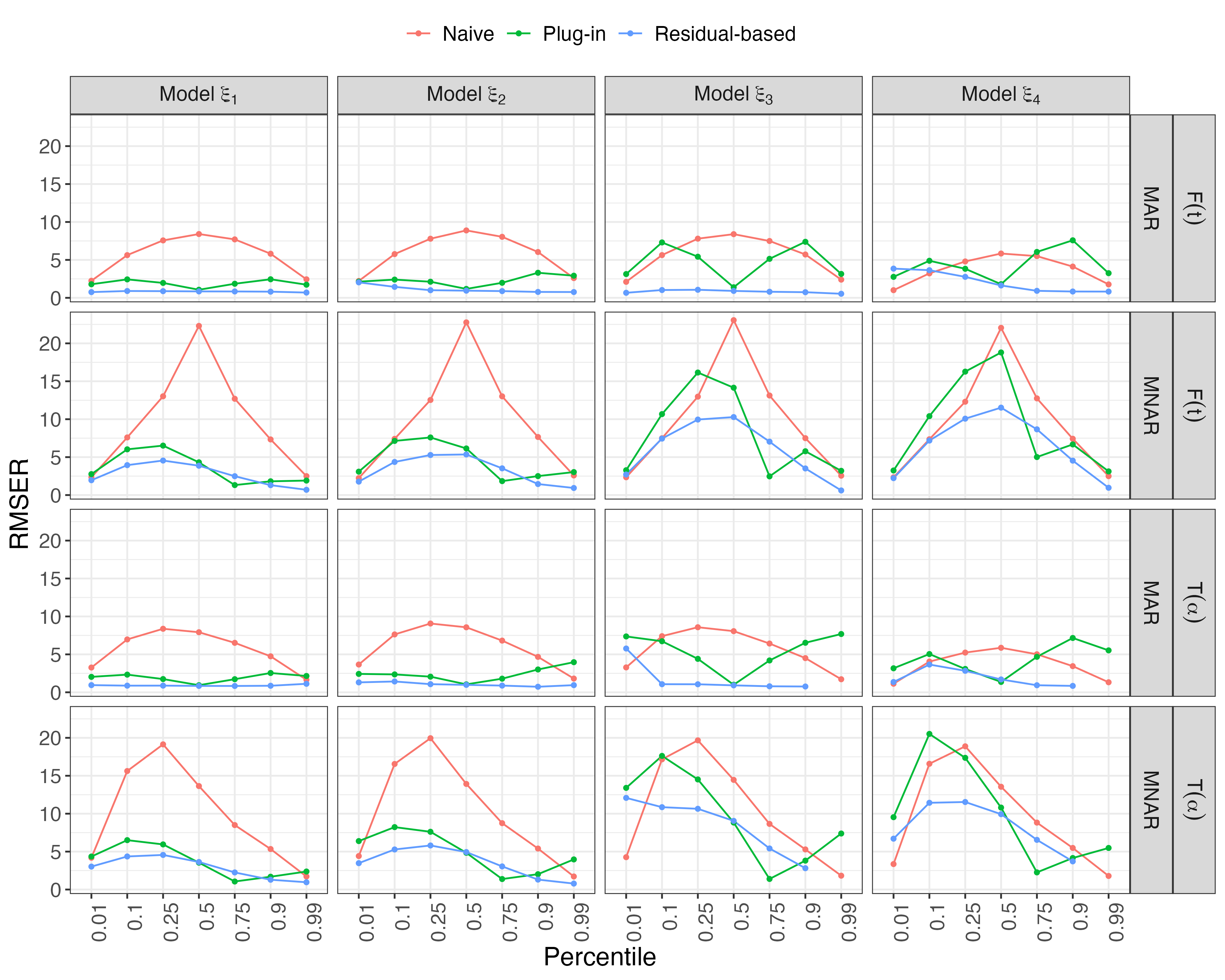}
   \caption{RMSER for the na\"{i}ve, plug-in, and residual-based CDF and quantile estimators when $n_\text{B} = n_\text{A}$. \label{fig:rmser_nA}}
 \end{figure}

\begin{figure}[H]
 \centering
  \includegraphics[width=0.9\textwidth,height=0.9\textheight,keepaspectratio]{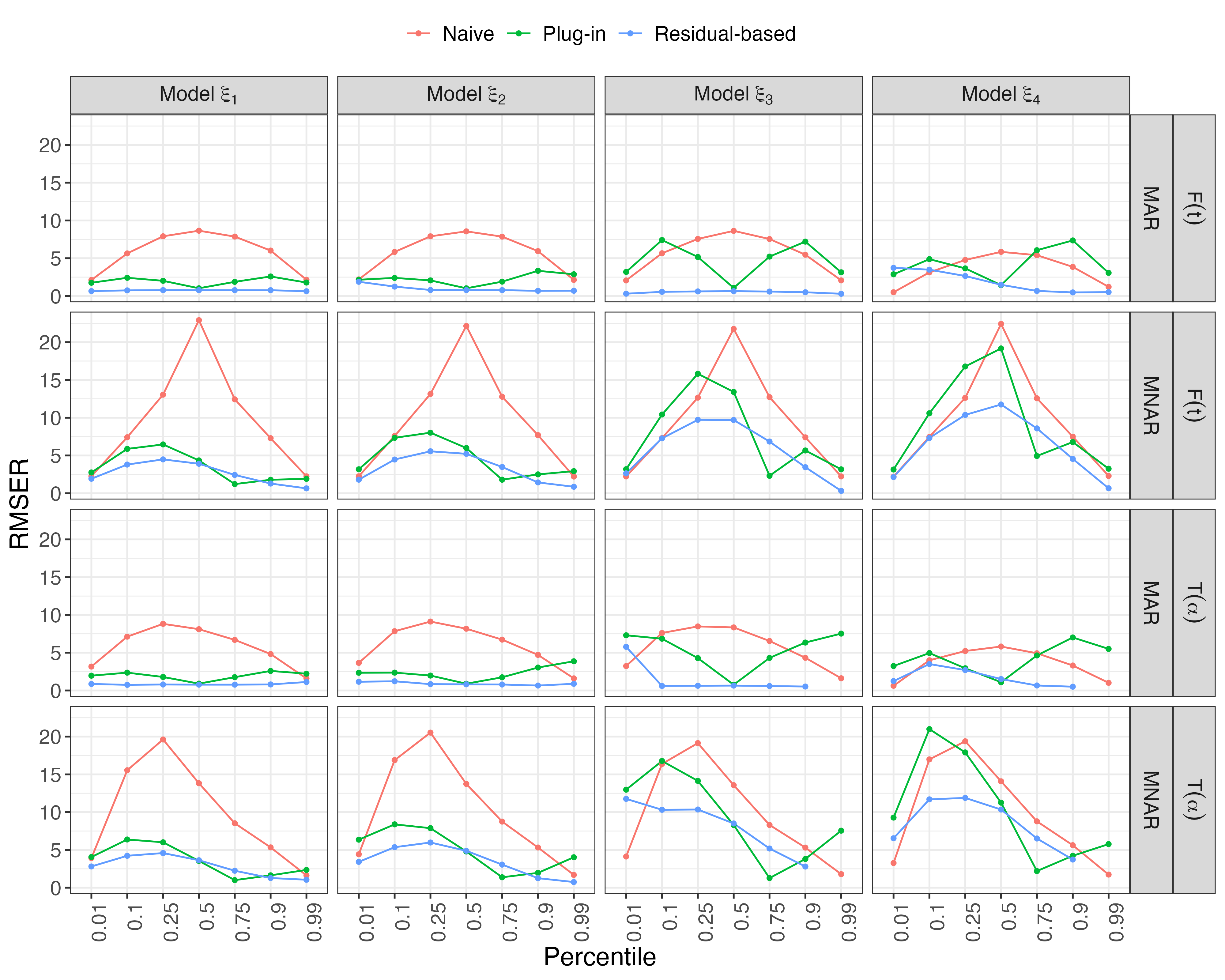}
   \caption{RMSER for the na\"{i}ve, plug-in, and residual-based CDF and quantile estimators when $n_\text{B} = 10n_\text{A}$. \label{fig:rmser_10nA}}
 \end{figure}

\begin{figure}[H]
  \centering
  \includegraphics[width=0.9\textwidth,height=0.9\textheight,keepaspectratio]{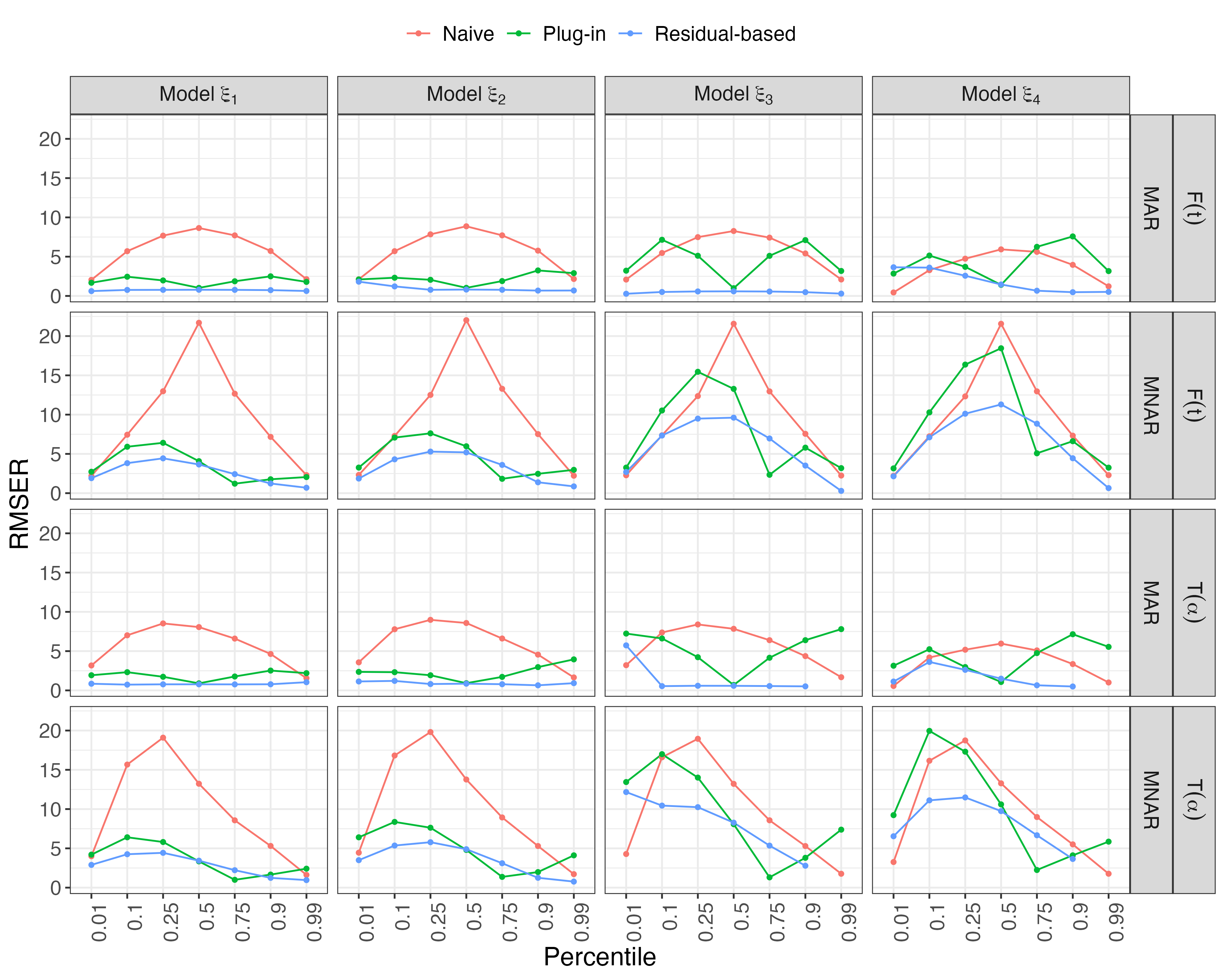}
   \caption{RMSER for the na\"{i}ve, plug-in, and residual-based CDF and quantile estimators when $n_\text{B} = 20n_\text{A}$. \label{fig:rmser_20nA}}
 \end{figure}

The simulation results for CDF and quantile estimation are presented in Figures \ref{fig:rmser_nA}--\ref{fig:rmser_20nA}. 

\subsection{CDF Estimation Results}
\label{sim:cdf}

Under MAR, the residual-based estimator $\widehat{F}_\text{R}(t; \boldsymbol{\widehat{\beta}})$ consistently outperformed both $\widehat{F}_\text{P}(t; \boldsymbol{\widehat{\beta}})$ and $\widehat{F}_\text{B}(t)$, exhibiting substantially lower RMSER values across nearly all models, percentiles, and $n_\text{B}$ combinations. The efficiency gains from $\widehat{F}_\text{R}(t; \boldsymbol{\widehat{\beta}})$ were most pronounced under models $\xi_{1}$–$\xi_{3}$. Notable exceptions occurred at the lower tail—particularly the 1st and 10th percentiles—under the highly nonlinear model $\xi_{4}$, where $\widehat{F}_\text{R}(t; \boldsymbol{\widehat{\beta}})$ had slightly elevated RMSER values.

Under MNAR, $\widehat{F}_\text{R}(t; \boldsymbol{\widehat{\beta}})$ remained the best-performing estimator overall, though the RMSER increased markedly for all estimators compared to MAR, reflecting the bias introduced by violations of the ignorability assumption. The loss of efficiency was especially evident under models $\xi_{3}$ and $\xi_{4}$, where the sample selection mechanism distorted the residual distribution and reduced the representativeness of $\mathcal{B}$. Nonetheless, $\widehat{F}_\text{R}(t; \boldsymbol{\widehat{\beta}})$ retained a clear advantage over its competitors across most percentiles and sample sizes.

Across both MAR and MNAR settings, the na\"{i}ve estimator $\widehat{F}_\text{B}(t)$ was, as expected, unaffected by the model specification since it ignores the covariates information. Its RMSER patterns were largely flat across models and only modestly improved with larger $n_\text{B}$. In contrast, the plug-in estimator $\widehat{F}_\text{P}(t; \boldsymbol{\widehat{\beta}})$ was notably more sensitive to model misspecification, with performance deteriorating sharply under nonlinear models where $m(\boldsymbol{X}; \boldsymbol{\widehat{\beta}})$ failed to capture the mean structure. Despite this general trend, two consistent patterns emerged. First, under MAR, the plug-in estimator $\widehat{F}_\text{P}(t; \boldsymbol{\widehat{\beta}})$ performed performed competitively near the median, likely because model misspecification has less effect near central quantiles where residual variability is smaller. Second, under MNAR, $\widehat{F}_\text{P}(t; \boldsymbol{\widehat{\beta}})$ showed relatively lower RMSER values at the upper–middle percentile ($\alpha = 0.75$), though this advantage diminished elsewhere. This pattern may reflect partial compensation between the selection bias in $\mathcal{B}$ and the regression adjustment based on $\boldsymbol{X}$, producing more accurate CDF alignment in that region of the distribution.

Finally, while $\widehat{F}_\text{B}(t)$ consistently exhibited the highest RMSER overall, its error levels were relatively stable at the distributional extremes (1st and 99th percentiles), regardless of sample size or selection mechanism. This stability suggests that $\widehat{F}_\text{B}(t)$ may still provide acceptable tail estimates of $F_\text{N}(t)$ in scenarios where only $\mathcal{B}$ is available.

\subsection{Quantile Estimation Results}
\label{sim:quant}

Overall, the performance of the three quantile estimators mirrored that of their CDF counterparts, with two important exceptions. First, $\widehat{T}_\text{R}(0.99)$ was unavailable in nearly all runs under models $\xi_{3}$ and $\xi_{4}$, regardless of sample size or sample selection mechanism. This outcome is expected: when computing $\widehat{F}_\text{P}(t; \boldsymbol{\widehat{\beta}})$ or $\widehat{F}_\text{B}(t)$, sample minima and maxima necessarily yield $\widehat{F}(t) = 0$ or $1$, even if the CDF is poorly estimated. In contrast, $\widehat{F}_\text{R}(t; \boldsymbol{\widehat{\beta}})$ may not reach these limits if no $\boldsymbol{X}_{i} \in \mathcal{A}$ yields $\widehat{G}\!\left(R_{i}(t; \boldsymbol{\widehat{\beta}})\right)$ close to zero or one—particularly under model misspecification. In such cases, substituting $\widehat{T}_\text{R}(\alpha)$ with $\widehat{T}_\text{B}(\alpha)$ may offer a more stable estimate of $T_\text{N}(\alpha)$.

A second notable exception occurred under MNAR for model $\xi_{4}$, where relative performance varied substantially across percentiles. Specifically, $\widehat{T}_\text{B}(\alpha)$ performed best at the distributional extremes (1st and 99th percentiles), $\widehat{T}_\text{P}(\alpha)$ at the 75th percentile (consistent with earlier findings for the CDF), and $\widehat{T}_\text{R}(\alpha)$ elsewhere. Despite $\widehat{F}_\text{R}(t; \boldsymbol{\widehat{\beta}})$ generally outperforming the other CDF estimators, the results suggest that, particularly at the distributional extremes, $\widehat{T}_\text{R}(\alpha)$ may be more sensitive to the joint effects of model misspecification and the nonignorable selection mechanism of sample $\mathcal{B}$—an issue that warrants further investigation in future research.

\subsection{Variance Estimation}
\label{sec:var}

We conclude this section with results from a limited variance estimation study that examined the performance of variance estimators for $\widehat{F}_\text{R}(t; \boldsymbol{\widehat{\beta}})$ and $\widehat{T}_\text{R}(\alpha)$ under $\boldsymbol{n}_\text{B} = [n_\text{A},\, 20n_\text{A}]$ and models $\xi_{1}$ and $\xi_{3}$. The simulation settings were otherwise identical to those described previously.

We computed two types of variance estimators for both $\widehat{F}_\text{R}(t; \boldsymbol{\widehat{\beta}})$ and $\widehat{T}_\text{R}(\alpha)$. The first, denoted $\widehat{V}_\text{Asymp}$, is given by Theorem~\ref{theorem_4} for $\widehat{F}_\text{R}(t; \boldsymbol{\widehat{\beta}})$ and by Eq.~\eqref{var_q} for $\widehat{T}_\text{R}(\alpha)$. The second, denoted $\widehat{V}_\text{Boot}$, is an adaptation of \textcite{kim2021combining}’s bootstrap variance estimation procedure and is summarized as follows:

\begin{enumerate}
\item Generate $L = 1{,}500$ sets of replicate weights for sample $\mathcal{A}$ and bootstrap samples (with replacement, of size $n_\text{B}$) from sample $\mathcal{B}$.
\item For each replicate, compute $\widehat{F}^{l}_\text{R}(t)$, $\widehat{T}^{l}_\text{R}(\alpha)$, and $\widehat{F}^{l}_\text{R}(\widehat{T}^{l}_\text{R}(\alpha))$.
\item Using the replicate estimates and the original estimates from $(\mathcal{A}, \mathcal{B})$, calculate  
\begin{align*}
\widehat{V}_\text{Boot}\!\left[\widehat{F}_\text{R}(t; \boldsymbol{\widehat{\beta}})\right] 
&= \frac{1}{L}\sum_{l = 1}^{L}\!\left(\widehat{F}^{l}_\text{R}(t) - \widehat{F}_\text{R}(t; \boldsymbol{\widehat{\beta}})\right)^2, \\
\widehat{V}_\text{Boot}\!\left[\widehat{T}_\text{R}(\alpha)\right] 
&= \left(\frac{\widehat{T}^{*}_\text{UL} - \widehat{T}^{*}_\text{LL}}{2z_{\gamma/2}}\right)^2,
\end{align*}
where  
\begin{align}
\widehat{T}^{*}_\text{LL} &= \inf\!\left\{t : \widehat{F}_\text{R}(t; \boldsymbol{\widehat{\beta}}) \geq \alpha - z_{\gamma/2}\widehat{V}^{1/2}_\text{Boot}\!\left[\widehat{F}_\text{R}\!\left(\widehat{T}_\text{R}(\alpha)\right)\right]\right\}, \label{wood_1}\\
\widehat{T}^{*}_\text{UL} &= \inf\!\left\{t : \widehat{F}_\text{R}(t; \boldsymbol{\widehat{\beta}}) \geq \alpha + z_{\gamma/2}\widehat{V}^{1/2}_\text{Boot}\!\left[\widehat{F}_\text{R}\!\left(\widehat{T}_\text{R}(\alpha)\right)\right]\right\} \label{wood_2},
\end{align}
with 
\begin{align*}
\widehat{V}_\text{Boot}\!\left[\widehat{F}_\text{R}\!\left(\widehat{T}_\text{R}(\alpha)\right)\right] 
&= \frac{1}{L}\sum_{l = 1}^{L}\!\left(\widehat{F}^{l}_\text{R}\left(\widehat{T}^{l}_\text{R}(\alpha)\right) - \widehat{F}_\text{R}\!\left(\widehat{T}_\text{R}(\alpha)\right)\right)^2. 
\end{align*}
\end{enumerate}

We then compared both $\widehat{V}_\text{Asymp}$ and $\widehat{V}_\text{Boot}$ with the ``gold-standard’’ Monte Carlo variance, defined as  
\begin{align*}
\widehat{V}_\text{MC} = \frac{1}{n_{\text{sim}} - 1}\sum_{m = 1}^{n_\text{sim}}\!\left(\hat{\theta}_{m} - \frac{1}{n_\text{sim}}\sum_{m = 1}^{n_\text{sim}}\!\hat{\theta}_{m}\right)^2,
\end{align*}
where $\hat{\theta}$ denotes either $\widehat{F}_\text{R}(t; \boldsymbol{\widehat{\beta}})$ or $\widehat{T}_\text{R}(\alpha)$.

After $n_\text{sim} = 1{,}500$ independent simulation runs, we computed the following performance statistics at the 1st, 10th, 25th, 50th, 75th, 90th, and 99th percentiles:

\begin{itemize}
\item \textbf{\% CR (Coverage Rate):} The percentage of instances in which the true population parameter was contained in the nominal 90\% $Z$-confidence interval.  
For $\widehat{F}_\text{R}(t; \boldsymbol{\widehat{\beta}})$, this was computed as  
\[
\widehat{F}_\text{R}(t; \boldsymbol{\widehat{\beta}}) \pm 1.645 \times \widehat{V}^{1/2}_\text{Asymp}\!\left[\widehat{F}_\text{R}(t; \boldsymbol{\widehat{\beta}})\right]
\]  
or  
\[
\widehat{F}_\text{R}(t; \boldsymbol{\widehat{\beta}}) \pm 1.645 \times \widehat{V}^{1/2}_\text{Boot}\!\left[\widehat{F}_\text{R}(t; \boldsymbol{\widehat{\beta}})\right],
\]  
and for $\widehat{T}_\text{R}(\alpha)$, using Eqs.~\eqref{q:LL}--\eqref{q:UL} or Eqs.~\eqref{wood_1}--\eqref{wood_2}.

\item \textbf{AL (Average Length):} The mean length of the confidence intervals across all $1{,}500$ simulations,  
\[
\text{AL} = \frac{1}{n_\text{sim}}\sum_{m = 1}^{n_\text{sim}}\!\left(\text{UL}_{m} - \text{LL}_{m}\right).
\]

\item \textbf{$\bar{\widehat{V}}$ (Average Estimated Variance):}  
\[
\bar{\widehat{V}} = \frac{1}{n_\text{sim}}\sum_{m = 1}^{n_\text{sim}}\!\widehat{V}_{m}.
\]

\item \textbf{\% RB (Percent Relative Bias):} The relative bias of the estimated variance with respect to $\widehat{V}_\text{MC}$,  
\[
\%\text{RB} = \frac{\bar{\widehat{V}} - \widehat{V}_\text{MC}}{\widehat{V}_\text{MC}} \times 100.
\]
\end{itemize}

\begin{figure}[ph]
    \centering
    \begin{subfigure}{\textwidth}
         \caption{Model $\xi_{1}$.}
        \includegraphics[width=\textwidth,height=\textheight,keepaspectratio]{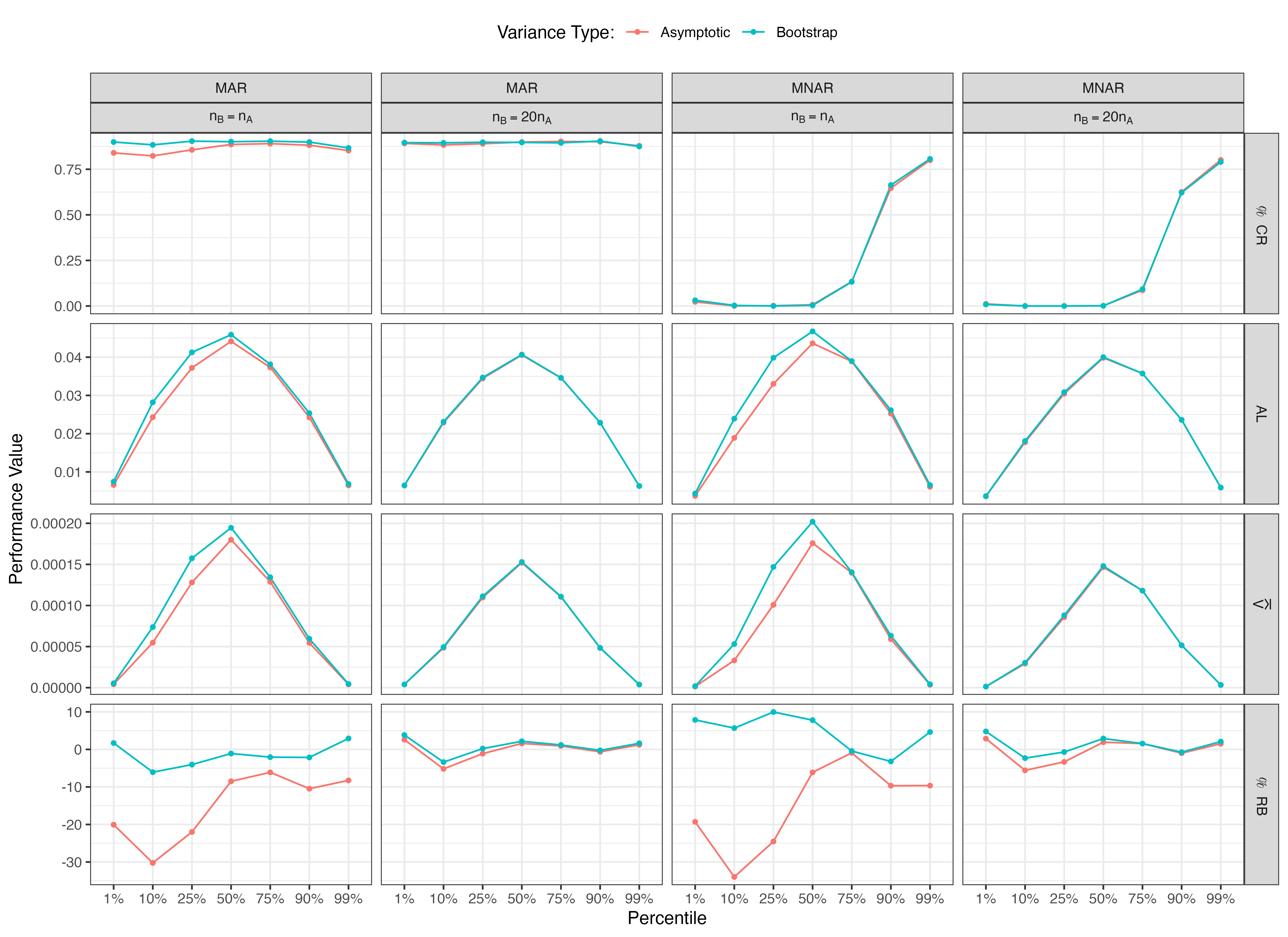}
    \end{subfigure}
    \vfill
    \begin{subfigure}{\textwidth}
   \caption{Model $\xi_{3}$.}
        \includegraphics[width=\textwidth,height=\textheight,keepaspectratio]{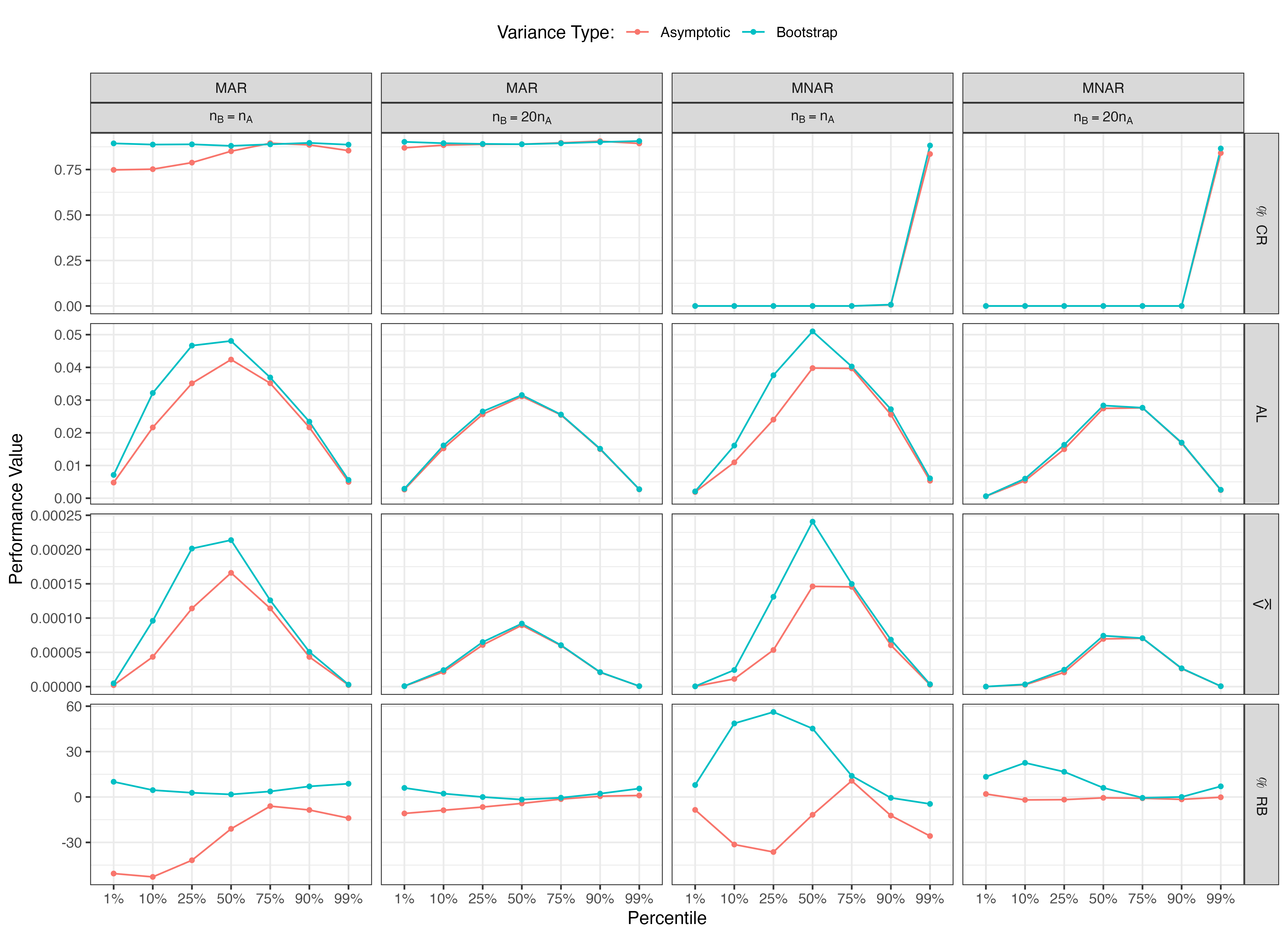}
    \end{subfigure}
    \caption{Performance statistics for the variance estimators of $\widehat{F}_\text{R}(t; \boldsymbol{\widehat{\beta}})$ under models $\xi_{1}$ and $\xi_{3}$, with $n_\text{A} = 1{\small,}000$.  \label{fig:ft}}
\end{figure}

\begin{figure}[ph]
    \centering
    \begin{subfigure}{1\textwidth}
        \caption{Model $\xi_{1}$.}
        \includegraphics[width=\textwidth,height=\textheight,keepaspectratio]{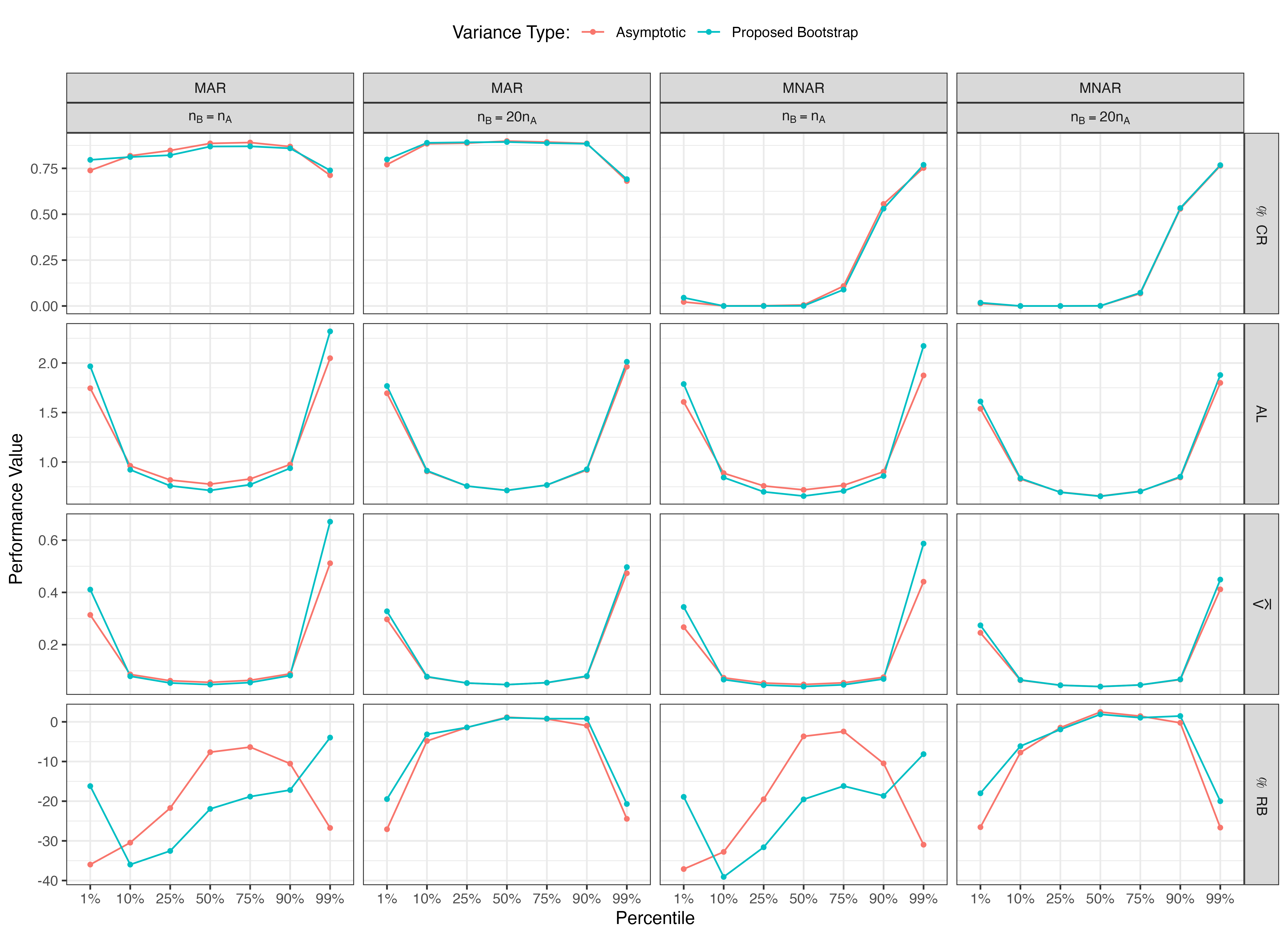}
    \end{subfigure}
    \vfill
    \begin{subfigure}{1\textwidth}
        \caption{Model $\xi_{3}$.}
        \includegraphics[width=\textwidth,height=\textheight,keepaspectratio]{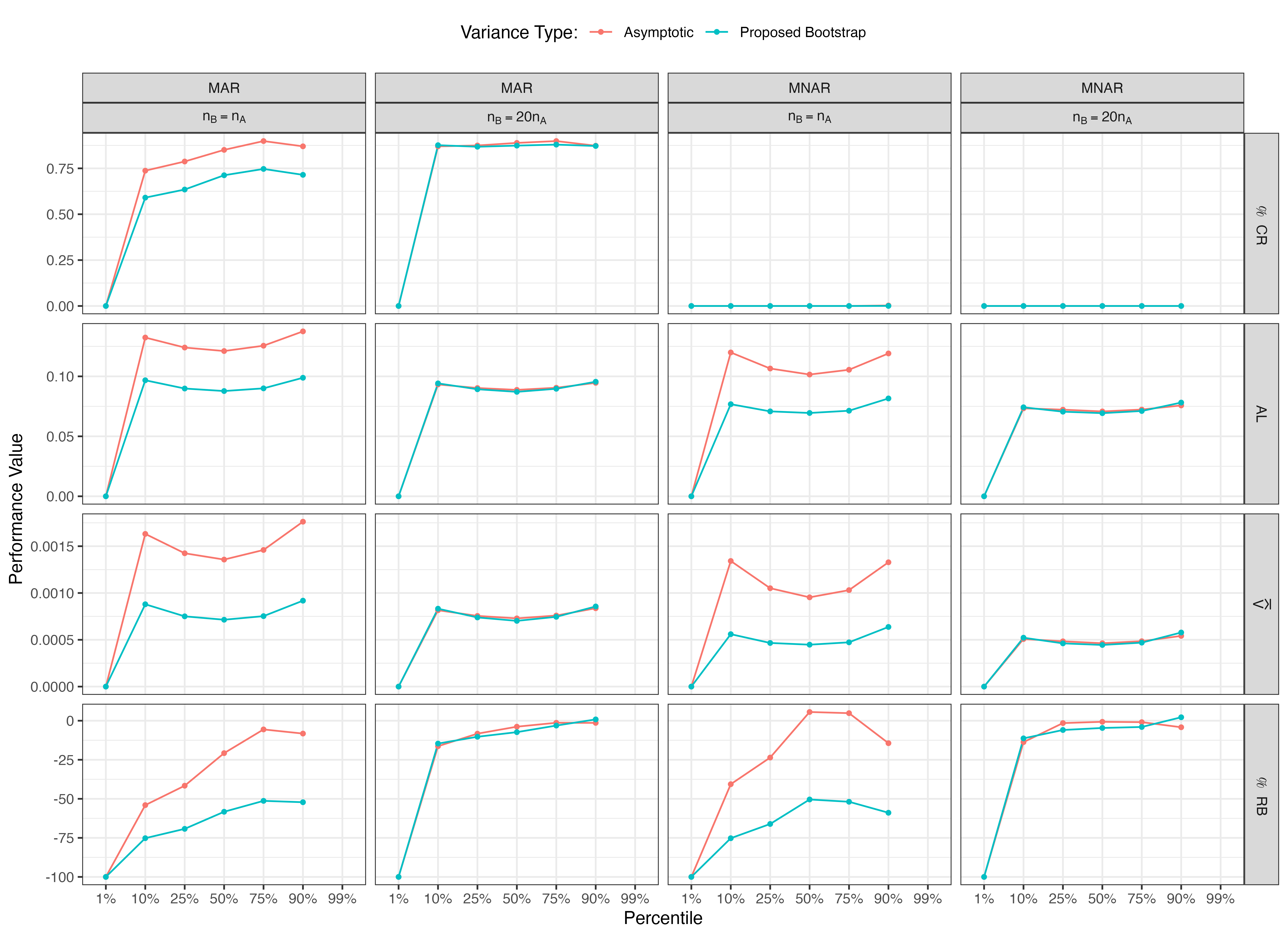}
    \end{subfigure}
    \caption{Performance statistics for the variance estimators of $\widehat{T}_\text{R}(\alpha)$ under models $\xi_{1}$ and $\xi_{3}$, with $n_\text{A} = 1{\small,}000$. \label{fig:ta}}
\end{figure}

The variance estimation results for $\widehat{F}_\text{R}(t; \boldsymbol{\widehat{\beta}})$ and $\widehat{T}_\text{R}(\alpha)$ are presented in Figures \ref{fig:ft} and \ref{fig:ta}, respectively.


\subsubsection{Variance Estimation Results for $\widehat{F}_{R}(t; \boldsymbol{\widehat{\beta}})$}

Starting with model~$\xi_{1}$, under MAR and $n_\text{B} = n_\text{A} = 1{\small,}000$, both variance estimators achieved coverage rates (CR) close to the nominal 90\%, though slightly lower for $\widehat{V}_\text{Asymp}$ due to its tendency to underestimate the true sampling variance, as reflected in the \%RB curves. As $n_\text{B}$ increased to $20{\small,}000$, the two estimators produced nearly identical results across all metrics. Under MNAR, both estimators exhibited substantially liberal coverage (CR well below 90\%), particularly for the lower percentiles (below the 75th), despite showing similar AL and \%RB patterns to those under MAR. This pattern suggests that the low coverage arises primarily from bias in $\widehat{F}_\text{R}(t; \boldsymbol{\widehat{\beta}})$ induced by the violation of the ignorability assumption, where 85\% of sample~$\mathcal{B}$ were constrained to fall above the median of the response variable.

For model~$\xi_{3}$, under MAR and $n_\text{B}= 1{\small,}000$, both variance estimators achieved good coverage rates, but $\widehat{V}_\text{Asymp}$ again tended to underestimate the sampling variance, with this bias more pronounced under model misspecification. The underestimation is evident from the lower coverage and higher \%RB, particularly at the lower percentiles. Increasing $n_\text{B}$ to $20{\small,}000$ nearly eliminated these discrepancies across percentiles. Under MNAR, neither estimator achieved satisfactory coverage due to the combined effects of nonignorable selection and model misspecification. Nevertheless, for smaller $n_\text{B}$, $\widehat{V}_\text{Boot}$ showed a clear overestimation bias, whereas $\widehat{V}_\text{Asymp}$ exhibited a milder underestimation bias. For larger $n_\text{B}$, both estimators exhibited substantial reductions in bias, with $\widehat{V}_\text{Asymp}$ displaying slightly greater stability in \%RB.

\subsubsection{Variance Estimation Results for $\widehat{T}_{R}(\alpha)$}

The variance estimation results for $\widehat{T}_\text{R}(\alpha)$ followed trends similar to those observed for $\widehat{F}_\text{R}(t; \boldsymbol{\widehat{\beta}})$, but with stronger sensitivity at the tails of the distribution. Under MAR and $n_\text{B}=n_\text{A}=1{\small,}000$, both $\widehat{V}_\text{Asymp}$ and $\widehat{V}_\text{Boot}$ tended to underestimate the true sampling variance across percentiles, resulting in coverage rates slightly below the nominal level. The underestimation bias was more pronounced at the lower percentiles (below $\alpha=0.50$) for both models~$\xi_{1}$ and $\xi_{3}$ and at the upper percentile ($\alpha=0.99$) for model~$\xi_{1}$. In general, $\widehat{V}_\text{Asymp}$ showed superior performance under model~$\xi_{3}$ in terms of coverage and \%RB. When $n_\text{B}$ increased to $20{\small,}000$, both estimators showed substantial reductions in bias and improved coverage, producing nearly identical CR, AL, and \%RB curves across percentiles. Nonetheless, substantial underestimation persisted at the extremes for both estimators under both models. Under MNAR, coverage deteriorated markedly—particularly at lower percentiles and under model~$\xi_3$—despite similar patterns of AL and \%RB, reflecting the effect of nonignorable selection bias carried over from $\widehat{F}_\text{R}(t; \boldsymbol{\widehat{\beta}})$.

These results reinforce that while both variance estimators perform well under MAR, large $n_\text{B}$, and correct model specification, their reliability diminishes under model misspecification, nonignorable selection, or both—particularly at the extreme quantiles. Because both estimators are derived from the \textcite{woodruff1952confidence} large-sample approximation, their performance tends to deteriorate when $n_\text{B}$ is small.

\section{Real Data Example}
\label{sec:realdata}

\noindent In this section, we illustrate the performance of $\widehat{F}_\text{R}(t; \boldsymbol{\widehat{\beta}})$ and $\widehat{T}_\text{R}(\alpha)$ using data from the National Health and Nutrition Examination Survey (NHANES) published by \textcite{CDC_NHANES_2015_2020}. NHANES is a multistage stratified random sample designed to obtain and disseminate detailed information on the health and nutritional status of a nationally representative sample of non-institutionalized individuals in the United States. Our goal is to estimate the population distribution function of total cholesterol (TCHL, $Y$, in mg/dL) among U.S. adults using biological sex ($X_1$), age ($X_2$), glycohemoglobin ($X_3$, in \%), triglycerides ($X_4$, in mg/dL), direct high-density lipoprotein cholesterol (HDL, $X_{5}$, in mg/dL), body mass index (BMI, $X_{6}$, in kg/$m^2$), and pulse ($X_{7}$) as covariates. We used the subsample of individuals with complete data from the 2015--2016 NHANES cohort as the probability sample ($n_\text{A} = 2{,}474$) and those from the 2017--2020 cohort as the nonprobability sample ($n_\text{B} = 3{,}770$). Notably, observations in $\mathcal{B}$ collected during the 2019--2020 round were deemed by the CDC to be nonrepresentative of the target population due to interruption in data collection during the COVID-19 pandemic and were combined with the previous cohort (2017-2018) to ensure respondent confidentiality.

As in Section~\ref{sim:cdf}, $\widehat{F}_\text{R}(t; \boldsymbol{\widehat{\beta}})$ and $\widehat{T}_\text{R}(\alpha)$ were constructed using a standard multivariable linear regression model fitted on sample $\mathcal{B}$ via \texttt{R}'s \texttt{lm()} function. We then computed the percent absolute relative bias of each estimator relative to its Horvitz–Thompson (HT) equivalent calculated using the 2015-2016 NHANES dataset (sample $\mathcal{A}$):
\begin{align}
\text{\%ARB}\!\left(\widehat{F}_\text{R}(t; \boldsymbol{\widehat{\beta}})\right) &= \left|\frac{\widehat{F}_\text{R}(t; \boldsymbol{\widehat{\beta}}) - \widehat{F}_{\pi}(t)}{\widehat{F}_{\pi}(t)}\right| \times 100, \\
\text{\%ARB}\!\left(\widehat{T}_\text{R}(\alpha)\right) &= \left|\frac{\widehat{T}_\text{R}(\alpha) - \widehat{T}_{\pi}(\alpha)}{\widehat{T}_{\pi}(\alpha)}\right| \times 100,
\end{align}
along with bootstrap estimates of their sampling variance. The use of the asymptotic variance estimator in Theorem \ref{theorem_4} was not considered here since the second-order inclusion probabilities were unavailable. 

\begin{figure}[h!]
\centering
\caption{Regression diagnostic plots: (a) residuals against the fitted values, (b) Normal quantile-quantile plot, (c) scale-location plot, and (d) standardized residuals against leverage values. \label{fig:diagn}}
\includegraphics[scale = .50]{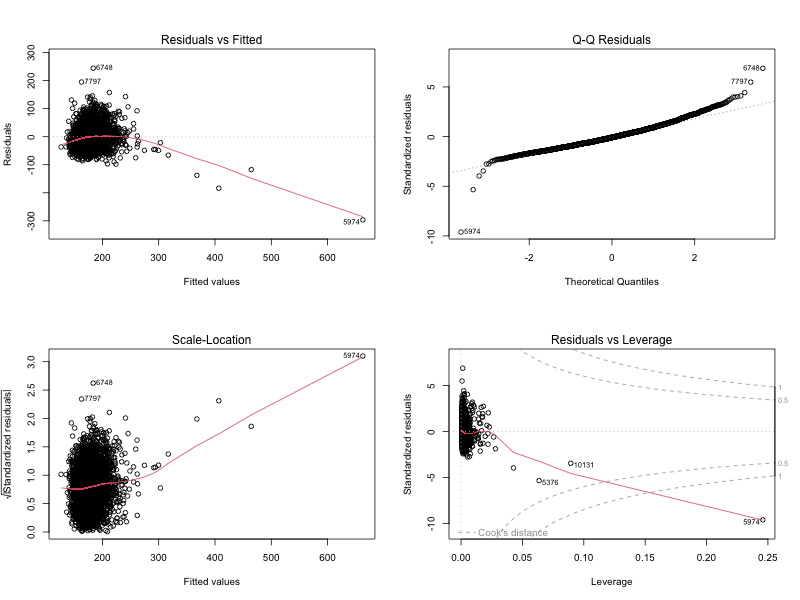}
\end{figure}

The results are presented in Tables~\ref{realdat}--\ref{boot_T}. To address uncertainty regarding the true regression function and the selection mechanism of sample $\mathcal{B}$, we also report regression coefficients (Table~\ref{tab:coeff}) and diagnostic plots (Figure~\ref{fig:diagn}). 

We begin by discussing some characteristics of the regression model fitted on $\mathcal{B}$. From Table~\ref{tab:coeff}, six out of eight (75\%) of the survey-weighted coefficients calculated from sample $\mathcal{A}$ using \textcite{thomlum}'s \texttt{svyglm()} function in \texttt{R} were captured within their respective 90\% confidence intervals constructed from the regression model fitted on the nonprobability sample $\mathcal{B}$. The coefficient for pulse ($\hat{\beta}_{7}$) had a large $p$-value ($0.67$), suggesting that, conditional on the remaining covariates, pulse was not a significant predictor of TCHL. Similarly, the coefficient for biological sex indicated only a small mean difference in TCHL between females and males. However, these interpretations should be treated cautiously given the likely model misspecification, as suggested by the small $R^2=0.26$ and the strong nonlinearity, heteroscedasticity, and outliers visible in Figure~\ref{fig:diagn}.

\begin{table}[h!]
  \centering
  \caption{Regression results: (a) survey-weighted $\hat{\beta}$ from $\mathcal{A}$ ($\hat{\beta}_{\mathcal{A}}$), (b) na\"{i}ve coefficients from $\mathcal{B}$ ($\hat{\beta}_{\mathcal{B}}$), (c) respective SE, 90\% $t$-confidence intervals, $t$-statistics, and $p$-values from $\mathcal{B}$, and (d) indicators for $\mathrm{LL} \leq \hat{\beta}_{\mathcal{A}} \leq \mathrm{UL}$. \label{tab:coeff}}
       \begin{tabular}{ccccccccc}\toprule
          & $\mathcal{A}$ & $\mathcal{B}$ & SE    & LL    & UL    & $t^*$   & $\Pr\left(T \geq |t^*|\right)$ & $\hat{\beta}_{\mathcal{A}} \in \mathrm{CI}$ \\
          \midrule
        $\hat{\beta}_{0}$ & 86.76 & 93.35 & 5.62  & 84.10 & 102.60 & 16.60 & 0.00  & Y \\
    $\hat{\beta}_{1}$: Female & 0.33  & 2.12  & 1.24  & 0.08  & 4.15  & 1.71  & 0.09  & Y \\
    $\hat{\beta}_{2}$ & 0.28  & 0.25  & 0.03  & 0.19  & 0.30  & 7.79  & 0.00  & Y \\
    $\hat{\beta}_{3}$ & -1.62 & -1.03 & 0.57  & -1.97 & -0.09 & -1.80 & 0.07  & Y \\
    $\hat{\beta}_{4}$ & 0.27  & 0.20  & 0.01  & 0.18  & 0.21  & 28.31 & 0.00  & N \\
    $\hat{\beta}_{5}$ & 0.94  & 1.00  & 0.04  & 0.93  & 1.07  & 22.57 & 0.00  & Y \\
    $\hat{\beta}_{6}$ & 0.43  & 0.25  & 0.09  & 0.11  & 0.39  & 2.90  & 0.00  & N \\
    $\hat{\beta}_{7}$ & 0.04  & -0.02 & 0.05  & -0.11 & 0.06  & -0.43 & 0.67  & Y \\
    \bottomrule
    \end{tabular}%
  \label{tab:addlabel}%
\end{table}%

\begin{table}[h!]
  \centering
  \caption{Percent $\mathrm{ARB}$ of the na\"{i}ve, plug-in, and residual-based CDF estimators, as well as their respective quantile estimators, relative to the HT equivalents $\widehat{F}_{\pi}(t)$ and $\widehat{T}_{\pi}(\alpha)$, from the NHANSE data. \label{realdat}}
     \begin{tabular}{ccccccccc} \toprule
          &       &       & \multicolumn{3}{c}{$\mathrm{\%ARB}\big(\widehat{F}\big)$} & \multicolumn{3}{c}{$\mathrm{\%ARB}\big(\widehat{T}\big)$} \\
          \cmidrule(lr){4-9}
    $\alpha$   & $\widehat{F}_{\pi}(t)$ & $\widehat{T}_{\pi}(\alpha)$ & B     & P     & R     & B     & P     & R \\
    \midrule
      1\%   & 0.01  & 106.00 & 76.89 & 100.00 & \textit{\textbf{42.41}} & \textit{\textbf{5.66}} & 40.02 & 26.69 \\
    10\%  & 0.10  & 138.00 & 50.17 & 99.41 & \textit{\textbf{17.34}} & 5.80  & 15.77 & \textit{\textbf{0.37}} \\
    25\%  & 0.26  & 158.00 & 29.79 & 70.54 & \textit{\textbf{9.92}} & 5.06  & 7.07  & \textit{\textbf{2.04}} \\
    50\%  & 0.51  & 184.00 & 16.54 & 17.15 & \textit{\textbf{7.02}} & 4.89  & \textit{\textbf{2.42}} & \textit{\textbf{2.42}} \\
    75\%  & 0.75  & 212.00 & 7.49  & 21.60 & \textit{\textbf{4.54}} & 3.77  & 8.23  & \textit{\textbf{2.44}} \\
    90\%  & 0.90  & 244.00 & 3.56  & 9.87  & \textit{\textbf{2.71}} & 4.51  & 14.24 & \textit{\textbf{3.60}} \\
    99\%  & 0.99  & 295.00 & 0.13  & 0.76  & \textit{\textbf{0.05}} & 0.68  & 18.69 & \textit{\textbf{0.50}} \\
    \bottomrule
    \end{tabular}%
  \label{tab:addlabel}%
\end{table}%

\begin{table}[h!]
  \centering
  \caption{Performance of the bootstrap variance estimator of $\widehat{F}_\text{R}(t; \boldsymbol{\widehat{\beta}})$ from the NHANES data. \label{boot_F}}
    \scalebox{.93}{
    \begin{tabular}{ccccccc} \toprule
    $\alpha$ & $\widehat{F}_{\pi}(t)$ & $\widehat{F}_\text{R}(t; \boldsymbol{\widehat{\beta}})$ & LL    & UL    & $\widehat{V}_{\text{Boot}}\times 10^3$ & $\widehat{F}_{\pi}(t) \in \text{CI}$ \\
    \midrule
      1\%   & 0.01  & 0.01  & 0.01  & 0.02  & 0.02  & Y \\
    10\%  & 0.10  & 0.12  & 0.11  & 0.14  & 0.08  & N \\
    25\%  & 0.26  & 0.28  & 0.26  & 0.30  & 0.17  & N \\
    50\%  & 0.51  & 0.55  & 0.52  & 0.57  & 0.24  & N \\
    75\%  & 0.75  & 0.79  & 0.76  & 0.81  & 0.19  & N \\
    90\%  & 0.90  & 0.93  & 0.92  & 0.94  & 0.06  & N \\
    99\%  & 0.99  & 0.99  & 0.99  & 0.99  & 0.01  & Y \\
    \bottomrule
    \end{tabular}}%
  \label{tab:addlabel}%
\end{table}%

\begin{table}[h!]
  \centering
   \caption{Performance of the bootstrap variance estimator of $\widehat{T}_\text{R}(\alpha)$ from the NHANES data. \label{boot_T}}
    \scalebox{.90}{\begin{tabular}{ccccccc} \toprule
    $\alpha$ & $\widehat{T}_{\pi}(\alpha)$ & $\widehat{T}_\text{R}(\alpha)$ & LL    & UL    & $\widehat{V}_{\text{Boot}}$ & $\widehat{T}_{\pi}(\alpha) \in \text{CI}$ \\
    \midrule
        1\%   & 106.00 & 134.29 & 134.29 & 134.29 & 0.00  & N \\
    10\%  & 138.00 & 137.49 & 134.29 & 138.50 & 1.64  & Y \\
    25\%  & 158.00 & 154.77 & 154.66 & 154.95 & 0.01  & N \\
    50\%  & 184.00 & 179.55 & 179.47 & 179.58 & 0.00  & N \\
    75\%  & 212.00 & 206.82 & 206.82 & 206.92 & 0.00  & N \\
    90\%  & 244.00 & 235.22 & 234.10 & 235.69 & 0.23  & N \\
    99\%  & 295.00 & 296.49 & 284.88 & 344.97 & 333.66 & Y \\
    \bottomrule
    \end{tabular}}%
  \label{tab:addlabel}%
\end{table}%

Despite the likely model misspecification, the residual-based estimator $\widehat{F}_\text{R}(t; \boldsymbol{\widehat{\beta}})$ exhibited substantially smaller relative bias than both the na\"{i}ve $\widehat{F}_\text{B}(t)$ and plug-in $\widehat{F}_\text{P}(t; \boldsymbol{\widehat{\beta}})$ estimators across nearly all percentiles, as seen from Table~\ref{realdat}. For the quantile estimators, $\widehat{T}_\text{R}(\alpha)$ performed best overall, though it was slightly outperformed by $\widehat{T}_\text{B}(\alpha)$ at the 1st percentile, consistent with the estimator’s known sensitivity to distributional extremes (see Section~\ref{sim:quant}). These findings are broadly consistent with our simulation results, reinforcing the robustness of $\widehat{F}_\text{R}(t; \boldsymbol{\widehat{\beta}})$ relative to both the na\"{i}ve and plug-in alternatives, even under model misspecification.

We conclude with the results for the bootstrap variance estimators under the NHANES data. For $\widehat{F}_\text{R}(t; \boldsymbol{\widehat{\beta}})$ (Table~\ref{boot_F}), the Horvitz–Thompson estimator $\widehat{F}_{\pi}(t)$ calculated from sample $\mathcal{A}$ was excluded from the bootstrap confidence intervals at most percentiles, likely due to the very tight intervals resulting from underestimation of the sampling variance. Similarly, for $\widehat{T}_\text{R}(\alpha)$ (Table~\ref{boot_T}), $\widehat{T}_{\pi}(\alpha)$ was excluded at the 1st, 25th, 50th, 75th, and 99th percentiles. These intervals, however, are intended to quantify uncertainty around $F_\text{N}(t)$ and $T_\text{N}(\alpha)$ rather than $\widehat{F}_{\pi}(t)$ or $\widehat{T}_{\pi}(\alpha)$, and should therefore be interpreted with caution.

\section{Discussion}
\label{sec:conclude}

\noindent This paper addressed the problem of estimating the finite population distribution function ($F_\text{N}(t)$) and quantiles under monotone missing data from combined probability and nonprobability samples. We proposed a residual-based estimator of $F_\text{N}(t)$ and established its large-sample properties. Under ignorability of the selection mechanism for the nonprobability sample and standard regularity conditions, the proposed estimator was shown to be asymptotically unbiased for $F_\text{N}(t)$. Simulation results demonstrated that the residual-based estimator consistently outperformed both the naïve empirical estimator (based solely on the nonprobability sample) and the pure mass imputation (plug-in) estimator across a range of models, percentiles, and sample size configurations. Under ignorability, it achieved substantially lower MSE values than its competitors, particularly under well-specified or moderately misspecified outcome models. When the ignorability assumption was violated, the performance of all estimators declined, yet the residual-based estimator maintained a clear advantage across most scenarios. The naïve estimator exhibited stable but inefficient behavior, whereas the plug-in estimator was more sensitive to model misspecification, performing competitively only near central quantiles. For quantile estimation, the overall trends mirrored those observed for the CDF estimators. However, the residual-based estimator was occasionally unavailable at extreme upper quantiles under highly nonlinear models, where the naïve estimator provided a more stable alternative. Under nonignorability, the relative performance varied across percentiles, with the naïve estimator performing better at the distributional extremes and the residual-based estimator retaining its advantage near central quantiles. These results suggest that although the residual-based quantile estimator generally offers clear efficiency gains, its performance at the extremes may be limited by the combined effects of model misspecification and nonignorable selection. We also proposed and evaluated two variance estimators for the residual-based CDF and quantile estimators: a linearization-based estimator derived from the asymptotic variance expression and a replication-based estimator using the bootstrap. Both achieved satisfactory coverage and relative bias under ignorability of the nonprobability sampling mechanism, large $n_\text{B}$, and correct model specification. However, their performance declined under small sample sizes, model misspecification, or nonignorable selection. Both estimators generally underestimated the true variance, with the asymptotic estimator exhibiting more pronounced underestimation, while the bootstrap estimator occasionally overestimated it when $n_\text{B}$ was small. Our future research will extend this framework to nonparametric regression settings, which can mitigate model dependence and improve robustness to misspecification, particularly in the tails of the distribution. 

\section*{Declarations}

\noindent {\bf Conflict of Interest}

\noindent No potential Conflict of interest was reported by the authors.

\noindent {\bf Data Availability Statement} 

\noindent The data supporting the findings of this study are publicly available from the U.S. Centers for Disease Control and Prevention (CDC) -- National Center for Health Statistics (NCHS) at \url{https://wwwn.cdc.gov/nchs/nhanes/Default.aspx}.

\printbibliography
\clearpage

\appendix

\renewcommand\thefigure{\thesection\arabic{figure}}

\section*{Appendix: Proofs}\label{append}

\renewcommand{\theequation}{A\arabic{equation}}
\setcounter{equation}{0} 

\begin{proof}[Proof of Lemma \ref{lemma_1}] 
Define  $\boldsymbol{\beta}^{*}_{s} \coloneqq \boldsymbol{\beta}^{*} + n^{-1/2}_\text{B}s$ for $s \in C_{s}$,  and consider
\[
Q_{n_\text{B}}(s)
:= \widehat{F}_\text{R}(t; \boldsymbol{\beta}^{*}_{s}) -\widehat{F}_\text{R}(t; \boldsymbol{\beta}^{*})
- F(t; \boldsymbol{\beta}^{*}_{s}) + F(t; \boldsymbol{\beta}^{*}).
\] By the triangle inequality,
\begin{align*}
|Q_{n_\text{B}}(s)|
&\leq \Big|\widehat{F}_\text{R}(t; \boldsymbol{\beta}^{*}_{s}) -\widehat{F}_\text{R}(t; \boldsymbol{\beta}^{*})
- F_\text{N}(t; \boldsymbol{\beta}^{*}_{s}) + F_\text{N}(t; \boldsymbol{\beta}^{*})\Big|  \\
& + \Big|F_\text{N}(t;\boldsymbol{\beta}^{*}_{s}) - F_\text{N}(t; \boldsymbol{\beta}^{*})
- F(t;\boldsymbol{\beta}^{*}_{s})+F(t; \boldsymbol{\beta}^{*})\Big| \label{p_a2} \\
&=: L_{1} + L_{2},
\end{align*}
where $L_2$ converges in probability to 0 by Assumption \ref{a6}.

It remains to show $L_{1}\xrightarrow{P}0$. For any  $0< \zeta < \frac{1}{2}$, partition the compact set $C_{s}$ as 
\[
C_{s} = \bigcup_{k = 1}^{N^{\zeta}} C_{s_{k}}.
\]
Now, for any set of $s_{k} \in C_{s_{k}}$, define
\begin{align*}
D_{1k} &= \Big| \{\widehat{F}_\text{R}(t; \boldsymbol{\beta}^{*}_{s_{k}}) - F_\text{N}(t;\boldsymbol{\beta}^{*}_{s_{k}})\}
- \{\widehat{F}_\text{R}(t; \boldsymbol{\beta}^{*}) - F_\text{N}(t;\boldsymbol{\beta}^{*})\}\Big| \\
 &= \left|\frac{1}{N}\sum_{u \in \mathcal{U}}{\frac{1}{N}\sum_{v \in \mathcal{U}}\Delta_{u,v}{\left\{\mathbbm{1}\left(\epsilon^{s^{*}_{k}}_{v} \leq \frac{t-m(\boldsymbol{X}_{u}; \boldsymbol{\beta}^{*}_{s_{k}})}{\nu(\boldsymbol{X}_{u})}\right) - \mathbbm{1}\left(\epsilon^{*}_{v} \leq \frac{t-m(\boldsymbol{X}_{u}; \boldsymbol{\beta}^{*})}{\nu(\boldsymbol{X}_{u})}\right) \right\}}}\right|
\end{align*}
and
\begin{align*}
D_{2k} &= \Big| \{\widehat{F}_\text{R}(t; \boldsymbol{\beta}^{*}_{s}) - F_\text{N}(t;\boldsymbol{\beta}^{*}_{s})\}
- \{\widehat{F}_\text{R}(t; \boldsymbol{\beta}^{*}_{s_{k}}) - F_\text{N}(t;\boldsymbol{\beta}^{*}_{s_{k}})\}\Big| \\
&=  \left|\frac{1}{N}\sum_{u \in \mathcal{U}}{\frac{1}{N}\sum_{v \in \mathcal{U}}\Delta_{u,v}{\left\{\mathbbm{1}\left(\epsilon^{s^{*}}_{v} \leq \frac{t-m(\boldsymbol{X}_{u}; \boldsymbol{\beta}^{*}_{s})}{\nu(\boldsymbol{X}_{u})}\right) - \mathbbm{1}\left(\epsilon^{s^{*}_{k}}_{v} \leq \frac{t-m\left(\boldsymbol{X}_{u}; \boldsymbol{\beta}^{*}_{s_{k}}\right)}{\nu(\boldsymbol{X}_{u})}\right) \right\}}}\right|,
\end{align*} 

where $$\boldsymbol{\beta}^{*}_{s_{k}} \coloneqq \boldsymbol{\beta}^{*} + n^{-1/2}_\text{B}s_{k},$$ $$\epsilon^{s^{*}_{k}}_{v} \coloneqq \frac{Y_{v}-m(\boldsymbol{X}_{v}; \boldsymbol{\beta}^{*}_{s_{k}})}{\nu(\boldsymbol{X}_{v})},$$  $$\Delta_{u,v} \coloneqq \left(\delta_{u}^{\mathcal{A}} \pi^{-1}_{u} \times \delta_{v}^{\mathcal{B}}f^{-1}_{n_\text{B}}\right) - 1,$$ and $$\delta_{s}^{\boldsymbol{S}} \coloneqq \begin{cases}
1 \ \mbox{if} \ s \in \boldsymbol{S} \\
0 \ \mbox{otherwise.}
\end{cases}
$$ Since $L_{1} \leq \max_{k}\{D_{1k}\}+\max_{k}\{D_{2k}\}$, it suffices to prove
\[
\max_{k}\{D_{1k}\} \xrightarrow{P}0
\quad\text{and}\quad
\max_{k}\{D_{2k}\}\xrightarrow{P}0.
\] For $D_{1k}$, note that \begin{align*} 
\Pr\left(\max_{k}\left\{D_{1k}\right\} \geq \varepsilon\right) &\leq  \sum_{k = 1}^{N^{\zeta}}{\Pr\left(D_{1k}\geq \varepsilon\right)}  \\
&\leq \frac{1}{\varepsilon}  \times\sum_{k = 1}^{N^{\zeta}}{\mathbb{E}\left(D_{1k}\right)}  \numberthis \label{markov_inequal} \\
&= \frac{1}{\varepsilon}  \times\sum_{k = 1}^{N^{\zeta}}{\mathbb{E}_{\mathscr{D}}\left[\mathbb{E}_{\xi}\left(D_{1k} \ | \ \delta^{\mathcal{A}}, \delta^{\mathcal{B}}, \boldsymbol{X}_\text{N}\right)\right]} \\
&\leq \frac{1}{\varepsilon}  \times \sum_{k = 1}^{N^{\zeta}}{{\max_{k}\left\{\mathbb{E}_{\mathscr{D}}\left[\mathbb{E}_{\xi}\left(D_{1k} \ | \ \delta^{\mathcal{A}}, \delta^{\mathcal{B}}, \boldsymbol{X}_\text{N}\right)\right]\right\}}}  \\
&\leq \frac{1}{\varepsilon} \times  N^{\zeta} \times \max_{k}\left\{\mathbb{E}_{\mathscr{D}}\left[\mathbb{E}_{\xi}\left(D_{1k} \ | \ \delta^{\mathcal{A}}, \delta^{\mathcal{B}}, \boldsymbol{X}_\text{N}\right)\right]\right\}
\end{align*} and
\begin{align*}
|D_{1k}| &\leq \frac{1}{N}\sum_{u \in \mathcal{U}}{\frac{1}{N}\sum_{v \in \mathcal{U}}{}}{\left|\Delta_{u,v}\right| \times \left|\mathbbm{1}\left(\epsilon^{s^{*}_{k}}_{v} \leq \frac{t-m(\boldsymbol{X}_{u}; \boldsymbol{\beta}^{*}_{s_{k}})}{\nu(\boldsymbol{X}_{u})}\right) - \mathbbm{1}\left(\epsilon^{*}_{v} \leq \frac{t-m(\boldsymbol{X}_{u}; \boldsymbol{\beta}^{*})}{\nu(\boldsymbol{X}_{u})}\right)\right|} \numberthis \label{triangle_inequal} \\
&\leq f^{-1}_{n_\text{B}} \times \frac{1}{N}\sum_{u \in \mathcal{U}}{\frac{\pi^{-1}_{u}}{N}\sum_{v \in \mathcal{U}}{}}{\left|\mathbbm{1}\left(\epsilon^{s^{*}_{k}}_{v} \leq \frac{t-m(\boldsymbol{X}_{u}; \boldsymbol{\beta}^{*}_{s_{k}})}{\nu(\boldsymbol{X}_{u})}\right) - \mathbbm{1}\left(\epsilon^{*}_{v} \leq \frac{t-m(\boldsymbol{X}_{u}; \boldsymbol{\beta}^{*})}{\nu(\boldsymbol{X}_{u})}\right)\right|} \numberthis \label{delta_max} \\
& \leq f^{-1}_{n_\text{B}} \times c \times \frac{1}{N}\sum_{u \in \mathcal{U}}{\frac{1}{N}\sum_{v \in \mathcal{U}}{}}{\left|\mathbbm{1}\left(\epsilon^{s^{*}_{k}}_{v} \leq \frac{t-m(\boldsymbol{X}_{u}; \boldsymbol{\beta}^{*}_{s_{k}})}{\nu(\boldsymbol{X}_{u})}\right) - \mathbbm{1}\left(\epsilon^{*}_{v} \leq \frac{t-m(\boldsymbol{X}_{u}; \boldsymbol{\beta}^{*})}{\nu(\boldsymbol{X}_{u})}\right)\right|} \numberthis \label{pi_max}
\end{align*} for some $c \in \mathbb{R}^{+}$, where Eq. \eqref{markov_inequal} follows from Markov's inequality (on the total expectation of $D_{1k}$), Eq. \eqref{triangle_inequal} from the triangle inequality,  Eq. \eqref{delta_max} from $\Delta_{u,v} \leq \pi^{-1}_{u}f^{-1}_{n_\text{B}} - 1 \leq \pi^{-1}_{u}f^{-1}_{n_\text{B}}$, and Eq. \eqref{pi_max} from Assumption \ref{a4A}. By the monotonicity property of expectations, 
$$|D_{1k}| \leq f^{-1}_{n_\text{B}} \times c \times \frac{1}{N}\sum_{u \in \mathcal{U}}{\frac{1}{N}\sum_{v \in \mathcal{U}}{}}{\left|\mathbbm{1}\left(\epsilon^{s^{*}_{k}}_{v} \leq \frac{t-m(\boldsymbol{X}_{u}; \boldsymbol{\beta}^{*}_{s_{k}})}{\nu(\boldsymbol{X}_{u})}\right) - \mathbbm{1}\left(\epsilon^{*}_{v} \leq \frac{t-m(\boldsymbol{X}_{u}; \boldsymbol{\beta}^{*})}{\nu(\boldsymbol{X}_{u})}\right)\right|} $$ implies
\begin{align*}
&\mathbb{E}_{\mathscr{D}}\left[\mathbb{E}_{\xi}\left(D_{1k} \ | \ \delta^{\mathcal{A}}, \delta^{\mathcal{B}}, \boldsymbol{X}_\text{N}\right)\right] \\ 
&\leq f^{-1}_{n_\text{B}} \times c \times \mathbb{E}_{\mathscr{D}}\left[\mathbb{E}_{\xi}\left(\frac{1}{N}\sum_{u \in \mathcal{U}}{\frac{1}{N}\sum_{v \in \mathcal{U}}{}}{\left|\mathbbm{1}\left(\epsilon^{s^{*}_{k}}_{v} \leq \frac{t-m(\boldsymbol{X}_{u}; \boldsymbol{\beta}^{*}_{s_{k}})}{\nu(\boldsymbol{X}_{u})}\right) - \mathbbm{1}\left(\epsilon^{*}_{v} \leq \frac{t-m(\boldsymbol{X}_{u}; \boldsymbol{\beta}^{*})}{\nu(\boldsymbol{X}_{u})}\right)\right|} \ \Bigg| \ \delta^{\mathcal{A}}, \delta^{\mathcal{B}}, \boldsymbol{X}_\text{N}\right)\right]  \\
&= f^{-1}_{n_\text{B}} \times c \times \mathbb{E}_{\xi}\left(\frac{1}{N}\sum_{u \in \mathcal{U}}{\frac{1}{N}\sum_{v \in \mathcal{U}}{}}{\left|\mathbbm{1}\left(\epsilon^{s^{*}_{k}}_{v} \leq \frac{t-m(\boldsymbol{X}_{u}; \boldsymbol{\beta}^{*}_{s_{k}})}{\nu(\boldsymbol{X}_{u})}\right) - \mathbbm{1}\left(\epsilon^{*}_{v} \leq \frac{t-m(\boldsymbol{X}_{u}; \boldsymbol{\beta}^{*})}{\nu(\boldsymbol{X}_{u})}\right)\right|} \ \Bigg| \ \boldsymbol{X}_\text{N}\right) \\
&= O\left(n^{-1/2}_\text{B}\right), \numberthis \label{ass_62}
\end{align*} where Eq. \eqref{ass_62} follows from Assumption \ref{a6}. Thus, 
\begin{align*} 
\Pr\left(\max_{k}\left\{D_{1k}\right\} \geq \varepsilon\right) &\leq\frac{1}{\varepsilon}  \times\sum_{k = 1}^{N^{\zeta}}{\mathbb{E}\left(D_{1k}\right)} \\
&\leq \frac{1}{\varepsilon} \times  N^{\zeta} \times \max_{k}\left\{\mathbb{E}_{\mathscr{D}}\left[\mathbb{E}_{\xi}\left(D_{1k} \ | \ \delta^{\mathcal{A}}, \delta^{\mathcal{B}}, \boldsymbol{X}_\text{N}\right)\right]\right\} \\
&= O\left(N^{\zeta} \times n^{-1/2}_\text{B}\right) \\
&=  O\left(f^{-1/2}_{n_\text{B}}\times N^{\zeta - 1/2}\right) \\
&=  O\left(N^{\zeta - 1/2}\right);
\end{align*} and since $\zeta < 1/2$, $\Pr\left(\max_{k}\left\{D_{1k}\right\} \geq \varepsilon\right)  \to 0$ as $N \to \infty$.

For $D_{2k}$, since $$\Pr\left(\max_{k}\left\{D_{2k}\right\}\geq \varepsilon \right) \leq \frac{1}{\varepsilon} \times N^{\zeta}\max_{k}\left\{\mathbb{E}\left(D_{2k}\right)\right\},$$ it suffices to show $N^{\zeta}\max_{k}\left\{\mathbb{E}\left(D_{2k}\right)\right\} \to 0$ as $N \to \infty$. From the triangle inequality, and noting that $\pi_{u}^{-1} \leq c < \infty$ for all $u \in \mathcal{U}$ (by Assumption \ref{a4A}), we observe 
\begin{align*}
D_{2k} &\leq f^{-1}_{n_\text{B}} \times c \times\frac{1}{N}\sum_{u \in \mathcal{U}}{\frac{1}{N}\sum_{v \in \mathcal{U}}{}}{\left|\mathbbm{1}\left(\epsilon^{s^{*}}_{v} \leq \frac{t-m(\boldsymbol{X}_{u}; \boldsymbol{\beta}^{*}_{s})}{\nu(\boldsymbol{X}_{u})}\right) - \mathbbm{1}\left(\epsilon^{*}_{v} \leq \frac{t-m(\boldsymbol{X}_{u}; \boldsymbol{\beta}^{*})}{\nu(\boldsymbol{X}_{u})}\right)\right|} \\
&+ f^{-1}_{n_\text{B}} \times c \times\frac{1}{N}\sum_{u \in \mathcal{U}}{\frac{1}{N}\sum_{v \in \mathcal{U}}{}}{\left|\mathbbm{1}\left(\epsilon^{*}_{v} \leq \frac{t-m(\boldsymbol{X}_{u}; \boldsymbol{\beta}^{*})}{\nu(\boldsymbol{X}_{u})}\right)-\mathbbm{1}\left(\epsilon^{s^{*}_{k}}_{v} \leq \frac{t-m(\boldsymbol{X}_{u}; \boldsymbol{\beta}^{*}_{s_{k}})}{\nu(\boldsymbol{X}_{u})}\right) \right|}.
\end{align*} Now we may use the monotonicity property of expectations (and Assumption \ref{a6}) to show
\begin{align*}
&\mathbb{E}_{\mathscr{D}}\left[\mathbb{E}_{\xi}\left(D_{2k} \ | \ \delta^{\mathcal{A}}, \delta^{\mathcal{B}}, \boldsymbol{X}_\text{N}\right)\right] \\
&\leq f^{-1}_{n_\text{B}} \times c \times\mathbb{E}_{\xi}\left(\frac{1}{N}\sum_{u \in \mathcal{U}}{\frac{1}{N}\sum_{v \in \mathcal{U}}{}}{\left|\mathbbm{1}\left(\epsilon^{s^{*}}_{v} \leq \frac{t-m(\boldsymbol{X}_{u}; \boldsymbol{\beta}^{*}_{s})}{\nu(\boldsymbol{X}_{u})}\right) - \mathbbm{1}\left(\epsilon^{*}_{v} \leq \frac{t-m(\boldsymbol{X}_{u}; \boldsymbol{\beta}^{*})}{\nu(\boldsymbol{X}_{u})}\right)\right|} \ \Bigg| \boldsymbol{X}_\text{N} \right) \\
&+ f^{-1}_{n_\text{B}} \times c \times \mathbb{E}_{\xi}\left(\frac{1}{N}\sum_{u \in \mathcal{U}}{\frac{1}{N}\sum_{v \in \mathcal{U}}{}}{\left|\mathbbm{1}\left(\epsilon^{*}_{v} \leq \frac{t-m(\boldsymbol{X}_{u}; \boldsymbol{\beta}^{*})}{\nu(\boldsymbol{X}_{u})}\right)-\mathbbm{1}\left(\epsilon^{s^{*}_{k}}_{v} \leq \frac{t-m(\boldsymbol{X}_{u}; \boldsymbol{\beta}^{*}_{s_{k}})}{\nu(\boldsymbol{X}_{u})}\right) \right|} \ \Bigg| \ \boldsymbol{X}_\text{N} \right) \\
&= O\left(n_\text{B}^{-1/2}\right);
\end{align*} and since
\begin{align*}
\Pr\left(\max_{k}\left\{D_{2k}\right\} \geq \varepsilon \right) = O\left(N^{\zeta- 1/2}\right),
\end{align*} we have that $\max_{k}\{D_{2k}\}\xrightarrow{P} 0$. Therefore, $Q_{n_{\text{B},s}} \xrightarrow{P} 0$.

Now, since $\boldsymbol{\widehat{\beta}} = \boldsymbol{\beta}^{*} + O_\text{P}(n^{-1/2}_\text{B})$ by Assumption \ref{a5A}, we can set $s =n^{1/2}_\text{B}\left(\boldsymbol{\widehat{\beta}} - \boldsymbol{\beta}^{*}\right)$ and obtain $$\widehat{F}_\text{R}(t; \boldsymbol{\widehat{\beta}}) -\widehat{F}_\text{R}(t; \boldsymbol{\beta}^{*})
- F(t; \boldsymbol{\widehat{\beta}}) + F(t; \boldsymbol{\beta}^{*}) = o_{P}(1).$$

The proof of $\left\{\widehat{F}_\text{R}(t; \boldsymbol{\widehat{\beta}}) -\widehat{F}_\text{R}(t; \boldsymbol{\beta}^{*})\right\} \xrightarrow{P} 0$ follows immediately because $\left\{F(t; \boldsymbol{\widehat{\beta}}) - F(t; \boldsymbol{\beta}^{*})\right\} \xrightarrow{P} 0$ by the continuity of $F(\cdot)$ and  $\boldsymbol{\widehat{\beta}} \xrightarrow{P} \boldsymbol{\beta}^{*}$.
\end{proof}

\begin{proof}[Proof of Theorem \ref{theorem_1}]
Recall that $\boldsymbol{X}_\text{N}$ denotes the $N \times p$ population covariate matrix and $\delta^{\mathcal{A}}, \delta^{\mathcal{B}}$ the sample indicators for $\mathcal{A}$ and $\mathcal{B}$. Using  Assumption \ref{E_FN} and the law of iterated expectation (with respect to the sampling design, $\mathscr{D}$, and the model, $\xi$), we obtain
\begin{align*}
\mathbb{E}\left(\widehat{F}_\text{R}(t; \boldsymbol{\beta}^{*}) - F_\text{N}(t)\right) &\coloneqq \mathbb{E}_{\mathscr{D}}\left[\mathbb{E}_\xi\left(\widehat{F}_\text{R}(t; \boldsymbol{\beta}^{*}) - F_\text{N}(t) \  \Big| \  \delta^{\mathcal{A}}, \delta^{\mathcal{B}}, \boldsymbol{X}_\text{N} \right)\right]\\ 
&= \mathbb{E}_{\mathscr{D}}\left[\frac{1}{N}\sum_{u \in \mathcal{U}}{\frac{\delta^{\mathcal{A}}_u}{\pi_{u}} {G\left(R_{u}(t, \boldsymbol{\beta^{*}})\right) - G\left(R_u(t;\boldsymbol{\beta})\right)}}\right]\\
&= \frac{1}{N}\sum_{u \in \mathcal{U}}\left\{{{G\left(R_{u}(t, \boldsymbol{\beta^{*}})\right) - G\left(R_u(t;\boldsymbol{\beta})\right)}}\right\}.
\end{align*} Hence, $$\mathrm{Bias}\!\left[\widehat{F}_\text{R}(t; \boldsymbol{\beta}^{*})\right] 
=  \frac{1}{N}\sum_{u \in \mathcal{U}}
\Big\{G\!\left(R_{u}(t, \boldsymbol{\beta^{*}})\right) - G\!\left(R_u(t;\boldsymbol{\beta})\right)\Big\}.$$ 

For the variance component, the law of total variance yields
\begin{align}
\label{var_1_part}
\mathrm{Var}\left[\widehat{F}_\text{R}(t; \boldsymbol{\beta}^{*})\right] &= \mathbb{E}_{\mathscr{D}}\left[\mathrm{Var}_{\xi}\left(\widehat{F}_\text{R}(t; \boldsymbol{\beta}^{*}) \ \big| \ \delta^{\mathcal{A}}, \delta^{\mathcal{B}}, \boldsymbol{X}_\text{N} \right)\right]  + \mathrm{Var}_{\mathscr{D}}\left[\mathbb{E}_{\xi}\left(\widehat{F}_\text{R}(t; \boldsymbol{\beta}^{*}) \ \big| \ \delta^{\mathcal{A}}, \delta^{\mathcal{B}}, \boldsymbol{X}_\text{N} \right)\right]. 
\end{align} 
Starting with the first term on the right hand side of \eqref{var_1_part}, and adopting the shorthand $R_{k} \coloneqq R_{k}(t; \boldsymbol{\beta}^{*})$, we note
\begin{align*}
\mathrm{Var}_{\xi}\left(\widehat{F}_\text{R}(t; \boldsymbol{\beta}^{*}) \ \big| \ \delta^{\mathcal{A}}, \delta^{\mathcal{B}}, \boldsymbol{X}_\text{N} \right) &= \mathbb{E}_{\xi}\left(\widehat{F}^{2}_\text{R}(t; \boldsymbol{\beta}^{*}) \ \big| \ \delta^{\mathcal{A}}, \delta^{\mathcal{B}}, \boldsymbol{X}_\text{N} \right) - \left\{\mathbb{E}_{\xi}\left(\widehat{F}_\text{R}(t; \boldsymbol{\beta}^{*}) \ \big| \ \delta^{\mathcal{A}}, \delta^{\mathcal{B}}, \boldsymbol{X}_\text{N} \right)\right\}^{2} \\
&= \frac{1}{N^2}\sum_{h \in \mathcal{A}}{\sum_{i \in \mathcal{A}}{\pi^{-1}_{h}\pi^{-1}_{i}\mathbb{E}_{\xi}\left(\widehat{G}(R_{h})\widehat{G}(R_{i})  \ \big| \ \delta^{\mathcal{A}}, \delta^{\mathcal{B}}, \boldsymbol{X}_\text{N} \right)}} \\
&- \left\{\mathbb{E}_{\xi}\left(\widehat{F}_\text{R}(t; \boldsymbol{\beta}^{*}) \ \big| \ \delta^{\mathcal{A}}, \delta^{\mathcal{B}}, \boldsymbol{X}_\text{N} \right)\right\}^{2} ,
\end{align*} where 
$$
\widehat{G}\left(R_{k}\right) = \frac{1}{n_\text{B}}\sum_{j \in \mathcal{B}}{\mathbbm{1}\left(\epsilon^{*}_{j} \leq R_{k}\right)}
$$ and 
\begin{align*}
\left\{\mathbb{E}_{\xi}\left(\widehat{F}_\text{R}(t; \boldsymbol{\beta}^{*}) \ \big| \ \delta^{\mathcal{A}}, \delta^{\mathcal{B}}, \boldsymbol{X}_\text{N} \right)\right\}^{2} &= \left\{\frac{1}{N}\sum_{i \in \mathcal{A}}{\pi^{-1}_{i}G\left(R_{i}\right)}\right\}^{2} \\
&= \frac{1}{N^2}\sum_{h \in \mathcal{A}}{\sum_{i \in \mathcal{A}}{\pi^{-1}_{h}\pi^{-1}_{i}G\left(R_{h}\right)G\left(R_{i}\right)}}.
\end{align*} To find $\mathbb{E}_{\xi}\left(\widehat{G}(R_{h})\widehat{G}(R_{i}) \ \big| \ \delta^{\mathcal{A}}, \delta^{\mathcal{B}}, \boldsymbol{X}_\text{N} \right)$, observe that $$\widehat{G}(R_{h})\widehat{G}(R_{i}) =\frac{1}{n^2_\text{B}}{\sum_{j \in \mathcal{B}}{\sum_{k\in\mathcal{B}}{\mathbbm{1}\left(\epsilon^{*}_{j} \leq R_{h}\right)\mathbbm{1}\left(\epsilon^{*}_{k} \leq R_{i}\right)}}}$$ can be written as 
\begin{align*}
\widehat{G}(R_{h})\widehat{G}(R_{i}) &= \frac{1}{n^{2}_\text{B}}\sum_{j \in \mathcal{B}}{\mathbbm{1}\left(\epsilon^{*}_{j} \leq R_{h}, \epsilon^{*}_{j} \leq R_{i}\right)} \\
&+ \frac{1}{n^{2}_\text{B}}\sum_{j \neq k}{\mathbbm{1}\left(\epsilon^{*}_{j} \leq R_{h}\right)\mathbbm{1}\left(\epsilon^{*}_{k} \leq R_{i}\right)}.
\end{align*}  Taking $\xi$-expectation to each component gives  
\begin{align*}
\mathbb{E}_{\xi}\left(\widehat{G}(R_{h})\widehat{G}(R_{i}) \ \big| \ \delta^{\mathcal{A}}, \delta^{\mathcal{B}}, \boldsymbol{X}_\text{N} \right) &= \frac{1}{n^{2}_\text{B}}\sum_{j \in \mathcal{B}}{\Pr\left(\epsilon^{*}_{j} \leq R_{h}, \epsilon^{*}_{j} \leq R_{i}\right)} + \frac{1}{n^{2}_\text{B}}\sum_{j \neq k}{\Pr\left(\epsilon^{*}_{j} \leq R_{h}\right)\Pr\left(\epsilon^{*}_{k} \leq R_{i}\right)} \\
&= \frac{1}{n_\text{B}}G\left(\min{\left\{R_{h}, R_{i}\right\}}\right)  + \frac{n_\text{B} - 1}{n_\text{B}}G(R_{h})G(R_{i})
\end{align*} since $\epsilon^{*}_{j}, \epsilon^{*}_{k}$ are i.i.d. by Assumption \ref{a3}. Thus,
\begin{align*}
\mathrm{Var}_{\xi}\left(\widehat{F}_\text{R}(t; \boldsymbol{\beta}^{*}) \ \big| \ \delta^{\mathcal{A}}, \delta^{\mathcal{B}}, \boldsymbol{X}_\text{N} \right) &= \frac{1}{N^2}\sum_{h \in \mathcal{A}}\sum_{i \in \mathcal{A}}{\pi^{-1}_{h}\pi^{-1}_{i}}\left\{\frac{1}{n_\text{B}}G\left(\min{\left\{R_{h}, R_{i}\right\}}\right)  +  \frac{n_\text{B} - 1}{n_\text{B}}G(R_{h})G(R_{i}) - G(R_{h})G(R_{i})\right\} \\
&= \frac{1}{n_\text{B}N^2}\sum_{h \in \mathcal{A}}\sum_{i \in \mathcal{A}}{\pi^{-1}_{h}\pi^{-1}_{i}}\left\{G\left(\min{\left\{R_{h}, R_{i}\right\}}\right)  -  G(R_{h})G(R_{i})\right\}, \numberthis \label{da_kill_p1}
\end{align*} and the design-based expectation of \eqref{da_kill_p1} yields
\begin{align*}
\mathbb{E}_{\mathscr{D}}\left[\mathrm{Var}_{\xi}\left(\widehat{F}_\text{R}(t; \boldsymbol{\beta}^{*}) \ \big| \ \delta^{\mathcal{A}}, \delta^{\mathcal{B}}, \boldsymbol{X}_\text{N} \right) \right] &= \frac{1}{n_\text{B}N^2}\sum_{u \in \mathcal{U}}\sum_{v \in \mathcal{U}}{\frac{\mathbb{E}_{\mathscr{D}}\left[\delta^{\mathcal{A}}_{u}\delta^{\mathcal{A}}_{v}\right]}{\pi_{u}\pi_{v}}}\left\{G\left(\min{\left\{R_{u}, R_{v}\right\}}\right)  -  G(R_{u})G(R_{v})\right\} \\
&= \frac{1}{n_\text{B}N^2}\sum_{u \in \mathcal{U}}\sum_{v \in \mathcal{U}}{\frac{\pi_{u,v}}{\pi_{u}\pi_{v}}}\left\{G\left(\min{\left\{R_{u}, R_{v}\right\}}\right)  -  G(R_{u})G(R_{v})\right\}. \numberthis \label{var_1}
\end{align*}
For the second term in \eqref{var_1_part}, since $$\mathbb{E}_{\xi}\left(\widehat{F}_\text{R}(t; \boldsymbol{\beta}^{*}) \ \big| \ \delta^{\mathcal{A}}, \delta^{\mathcal{B}}, \boldsymbol{X}_\text{N} \right)= \frac{1}{N}\sum_{i \in \mathcal{A}}{\pi^{-1}_{i}G\left(R_{i}\right)},$$ it immediately follows that 
\begin{align}
\mathrm{Var}_{\mathscr{D}}\left[\mathbb{E}_{\xi}\left(\widehat{F}_\text{R}(t; \boldsymbol{\beta}^{*}) \ \big| \ \delta^{\mathcal{A}}, \delta^{\mathcal{B}}, \boldsymbol{X}_\text{N} \right)\right] = \frac{1}{N^2}\sum_{u \in \mathcal{U}}{\sum_{v \in \mathcal{U}}{\left\{\frac{\pi_{u,v}}{\pi_{u}\pi_{v}}-1\right\}G\left(R_{u}\right)G\left(R_{v}\right)}}. \label{var_2}
\end{align} Thus, adding \eqref{var_1} and \eqref{var_2} together, we have 
\begin{align*}
\mathrm{Var}\left[\widehat{F}_\text{R}(t; \boldsymbol{\beta}^{*})\right] &= \frac{1}{n_\text{B}N^2}\sum_{u \in \mathcal{U}}\sum_{v \in \mathcal{U}}{\frac{\pi_{u,v}}{\pi_{u}\pi_{v}}}\left\{G\left(\min{\left\{R_{u}, R_{v}\right\}}\right)  -  G(R_{u})G(R_{v})\right\} \\ 
&+ \frac{1}{N^2}\sum_{u \in \mathcal{U}}{\sum_{v \in \mathcal{U}}{\left\{\frac{\pi_{u,v}}{\pi_{u}\pi_{v}}-1\right\}G\left(R_{u}\right)G\left(R_{v}\right)}}.
\end{align*}
\end{proof}

\begin{proof}[Proof of Theorem \ref{theorem_2}]
From Lemma \ref{lemma_1}, if Assumptions \ref{a1}-\ref{a6} hold, then \[
\widehat{F}_\text{R}(t; \boldsymbol{\widehat{\beta}}) - 
\widehat{F}_\text{R}(t; \boldsymbol{\beta}^{*}) \xrightarrow{P} 0.
\] Recall from Assumption \ref{a4A} that there exists a $c \in \mathbb{R}$ such that $\pi^{-1}_{u} \leq c$ for all $u \in \mathcal{U}$. Since \begin{align*}
\left|\widehat{F}_\text{R}(t; \boldsymbol{\widehat{\beta}}) - 
\widehat{F}_\text{R}(t; \boldsymbol{\beta^{*}})\right|
&= \left|\frac{1}{N}\sum_{i \in \mathcal{A}} \pi_i^{-1} 
   \Big\{\widehat{G}(R_i(t;\boldsymbol{\widehat{\beta}}))
          - \widehat{G}(R_i(t;\boldsymbol{\beta^{*}})) \Big\}\right|\\
&\leq \frac{1}{N}\sum_{i \in \mathcal{A}} \pi_i^{-1} 
   \left|\widehat{G}(R_i(t;\boldsymbol{\widehat{\beta}}))
          - \widehat{G}(R_i(t;\boldsymbol{\beta^{*}})) \right|\\
   &\leq c \frac{n_\text{A}}{N},
\end{align*} it follows that the random sequence $\widehat{F}_\text{R}(t; \boldsymbol{\widehat{\beta}}) - 
\widehat{F}_\text{R}(t; \boldsymbol{\beta^{*}})$ is bounded in probability. Hence, by the dominated convergence theorem, $$\mathbb{E}\!\left(\widehat{F}_\text{R}(t; \boldsymbol{\widehat{\beta}})\right)- \mathbb{E}\left(\widehat{F}_\text{R}(t; \boldsymbol{\beta^{*}})\right)  \to 0.$$ It now remains to show that \begin{align}
\mathrm{Var}\!\left(\widehat{F}_\text{R}(t; \boldsymbol{\widehat{\beta}})\right)
- \mathrm{Var}\!\left(\widehat{F}_\text{R}(t; \boldsymbol{\beta^{*}})\right) \to 0. \label{var_diffs} \end{align} Since the left hand side of \eqref{var_diffs} can be written as 
\begin{align*}
\mathrm{Var}\!\left(\widehat{F}_\text{R}(t; \boldsymbol{\widehat{\beta}})\right)
- \mathrm{Var}\!\left(\widehat{F}_\text{R}(t; \boldsymbol{\beta^{*}})\right) &= \left\{\mathbb{E}\!\left(\widehat{F}_\text{R}(t;\boldsymbol{\widehat{\beta}})\right)\right\}^2
      - \left\{\mathbb{E}\!\left(\widehat{F}_\text{R}(t;\boldsymbol{\beta}^{*}) \right)\right\}^2 \\
      &+ \mathbb{E}\!\left(\widehat{F}^{2}_\text{R}(t;\boldsymbol{\beta}^{*})\right) 
      - \mathbb{E}\!\left(\widehat{F}^{2}_\text{R}(t;\boldsymbol{\widehat{\beta}})\right) \\
      &= S_{1}  + S_{2},
\end{align*} it suffices to show $S_{1} \to 0$ and $S_{2} \to 0$. 

$S_{1}$, as a difference of squares, can be written as
\begin{align*}
 \left\{\mathbb{E}\!\left(\widehat{F}_\text{R}(t;\boldsymbol{\widehat{\beta}})\right)\right\}^2
      - \left\{\mathbb{E}\!\left(\widehat{F}_\text{R}(t;\boldsymbol{\beta}^{*}) \right)\right\}^2  &= \left\{\mathbb{E}\!\left(\widehat{F}_\text{R}(t;\boldsymbol{\widehat{\beta}})\right) - \mathbb{E}\!\left(\widehat{F}_\text{R}(t;\boldsymbol{\beta}^{*}) \right) \right\} \\
      &\times \left\{\mathbb{E}\!\left(\widehat{F}_\text{R}(t;\boldsymbol{\widehat{\beta}})\right) + \mathbb{E}\!\left(\widehat{F}_\text{R}(t;\boldsymbol{\beta}^{*}) \right) \right\} \numberthis \label{sum_exp},
\end{align*} where the first factor in \eqref{sum_exp} converges to 0 by the dominated convergence theorem. The second factor in \eqref{sum_exp} is bounded by a real constant, since
\begin{align*}
\mathbb{E}\!\left(\widehat{F}_\text{R}(t;\boldsymbol{\widehat{\beta}})\right) + \mathbb{E}\!\left(\widehat{F}_\text{R}(t;\boldsymbol{\beta}^{*}) \right) &= \mathbb{E}\left(\frac{1}{N}\sum_{i \in \mathcal{A}}{\pi^{-1}_{i}\left\{\widehat{G}(R_{h}(t; \boldsymbol{\widehat{\beta}})) + \widehat{G}(R_{h}(t; \boldsymbol{\beta^{*}})) \right\}}\right) \\
&\leq \mathbb{E}\left(2c \times \frac{n_\text{A}}{N}\right)
\end{align*} by Assumption \ref{a4A} and the monotonicity property of expectations.  Therefore, since  

\[
\mathbb{E}\!\left(\widehat{F}_\text{R}(t;\boldsymbol{\widehat{\beta}})\right) - \mathbb{E}\!\left(\widehat{F}_\text{R}(t;\boldsymbol{\beta}^{*}) \right) \to 0 \quad \text{and} \quad 
\mathbb{E}\!\left(\widehat{F}_\text{R}(t;\boldsymbol{\widehat{\beta}})\right) + \mathbb{E}\!\left(\widehat{F}_\text{R}(t;\boldsymbol{\beta}^{*}) \right) < \infty,
\] it follows immediately that $$ \left\{\mathbb{E}\!\left(\widehat{F}_\text{R}(t;\boldsymbol{\widehat{\beta}})\right)\right\}^2
      - \left\{\mathbb{E}\!\left(\widehat{F}_\text{R}(t;\boldsymbol{\beta}^{*}) \right)\right\}^2 \to 0.$$
To prove $S_{2} \to 0$ by the dominated convergence theorem, it suffices to show 
\[
\widehat{F}^{2}_\text{R}(t; \boldsymbol{\widehat{\beta}}) - \widehat{F}^{2}_\text{R}(t; \boldsymbol{\beta^{*}}) \xrightarrow{P} 0 \quad \text{and} \quad \widehat{F}^{2}_\text{R}(t; \boldsymbol{\widehat{\beta}}) - \widehat{F}^{2}_\text{R}(t; \boldsymbol{\beta^{*}}) < \infty.
\] Observe that $\widehat{F}^{2}_\text{R}(t; \boldsymbol{\widehat{\beta}}) - 
\widehat{F}^{2}_\text{R}(t; \boldsymbol{\beta^{*}})$ can be factored as 
\begin{align*}
\widehat{F}^{2}_\text{R}(t; \boldsymbol{\widehat{\beta}}) - 
\widehat{F}^{2}_\text{R}(t; \boldsymbol{\beta^{*}})
&= \left(\widehat{F}_\text{R}(t; \boldsymbol{\widehat{\beta}}) - 
       \widehat{F}_\text{R}(t; \boldsymbol{\beta^{*}})\right) 
   \left(\widehat{F}_\text{R}(t; \boldsymbol{\widehat{\beta}}) + 
       \widehat{F}_\text{R}(t; \boldsymbol{\beta^{*}})\right) \\
&=: A_{n_\text{B}} \times B_{n_\text{B}}.
\end{align*}
 $A_{n_\text{B}}$ converges in probability to 0 by Lemma \ref{lemma_1}, and 
\begin{align*}
B_{n_\text{B}} 
&= \frac{1}{N}\sum_{i \in A} \pi_i^{-1}
   \Big[\widehat{G}\big(R_i(t;\boldsymbol{\widehat{\beta}})\big) +
        \widehat{G}\big(R_i(t;\boldsymbol{\beta^{*}})\big)\Big] \\
&\leq 2c\frac{n_\text{A}}{N}
\end{align*}
is bounded in probability. Thus, by Lemma 2.3.1 of  \textcite{lehmann1999elements}, the product $A_{n_\text{B}} \times B_{n_\text{B}} $ converges in probability to 0. It only remains to note that 
\begin{align*}
\left|\widehat{F}^{2}_\text{R}(t; \boldsymbol{\widehat{\beta}}) - \widehat{F}^{2}_\text{R}(t; \boldsymbol{\beta^{*}})\right| &\leq \frac{1}{N^2}\sum_{h \in \mathcal{A}}{\sum_{i \in \mathcal{A}}{\pi^{-1}_{h}\pi^{-1}_{i}}\left|\widehat{G}(R_{h}(t; \boldsymbol{\widehat{\beta}}))\widehat{G}(R_{i}(t; \boldsymbol{\widehat{\beta}})) - \widehat{G}(R_{h}(t; \boldsymbol{\beta}^{*}))\widehat{G}\left(R_{i}(t; \boldsymbol{\beta}^{*})\right)\right|} \\
&\leq \left(c\frac{n_\text{A}}{N}\right)^2,
\end{align*} and thus $\widehat{F}^{2}_\text{R}(t; \boldsymbol{\widehat{\beta}}) - \widehat{F}^{2}_\text{R}(t; \boldsymbol{\beta^{*}})$ is bounded. This completes the proof. 
\end{proof}

\begin{proof}[Proof of Theorem \ref{theorem_3}]
By the law of iterated expectation, $\mathbb{E}\left(\widetilde{V}(t; \boldsymbol{\beta}^{*})\right)$ can be rewritten as
\begin{align}
\mathbb{E}\left(\widetilde{V}(t; \boldsymbol{\beta}^{*})\right) &= \mathbb{E}_{\mathscr{D}}\left[{\mathbb{E}_{\xi}}\left(\widetilde{V}_{1} \ \big| \ \delta^{\mathcal{A}}, \delta^{\mathcal{B}}, \boldsymbol{X}_\text{N} \right)\right] \label{unbiased_varest_1} \\
&+ \mathbb{E}_{\mathscr{D}}\left[{\mathbb{E}_{\xi}}\left(\widetilde{V}_{2} \ \big| \ \delta^{\mathcal{A}}, \delta^{\mathcal{B}}, \boldsymbol{X}_\text{N} \right)\right] \label{unbaised_varest_2}.
\end{align} 
We start by noting that the inner expectation of \eqref{unbiased_varest_1} can be expanded as follows: 
\begin{align*}
\mathbb{E}_{\xi}\left(\widetilde{V}_{1} \ \big| \ \delta^{\mathcal{A}}, \delta^{\mathcal{B}}, \boldsymbol{X}_\text{N} \right) &= \frac{1}{(n_\text{B} - 1)N^2}\sum_{h \in \mathcal{A}}{\sum_{i \in \mathcal{A}}{\pi^{-1}_{h,i}\left\{\frac{\pi_{h,i}}{\pi_{h}\pi_{i}}-1\right\}\mathbb{E}_{\xi}\left\{n_\text{B}\widehat{G}\left(R_{h}\right)\widehat{G}\left(R_{i}\right) \ \big| \ \delta^{\mathcal{A}}, \delta^{\mathcal{B}}, \boldsymbol{X}_\text{N}\right\}}} \\
&- \frac{1}{(n_\text{B} - 1)N^2}\sum_{h \in \mathcal{A}}{\sum_{i \in \mathcal{A}}{\pi^{-1}_{h,i}\left\{\frac{\pi_{h,i}}{\pi_{h}\pi_{i}}-1\right\}\mathbb{E}_{\xi}\left\{\widehat{G}\left(\min\left\{R_{h}, R_{i}\right\}\right)\ \big| \ \delta^{\mathcal{A}}, \delta^{\mathcal{B}}, \boldsymbol{X}_\text{N}\right\}}}.
\end{align*} By the proof of Theorem \ref{theorem_1}, 
\begin{align*}
\mathbb{E}_{\xi}\left[\widehat{G}(R_{k}) \ \big| \ \delta^{\mathcal{A}}, \delta^{\mathcal{B}}, \boldsymbol{X}_\text{N} \right] &= G\left(R_{k}\right) \\
\mathbb{E}_{\xi}\left[\widehat{G}(R_{k})\widehat{G}(R_{l}) \ \big| \ \delta^{\mathcal{A}}, \delta^{\mathcal{B}}, \boldsymbol{X}_\text{N} \right] &= \frac{1}{n_\text{B}}G\left(\min{\left\{R_{k}, R_{l}\right\}}\right)  + \frac{n_\text{B} - 1}{n_\text{B}}G(R_{k})G(R_{l}),
\end{align*} and thus $\mathbb{E}_{\xi}\left(\widetilde{V}_{1} \ \big| \ \delta^{\mathcal{A}}, \delta^{\mathcal{B}}, \boldsymbol{X}_\text{N} \right)$ can be simplified as
\begin{align*}
\mathbb{E}_{\xi}\left(\widetilde{V}_{1} \ \big| \ \delta^{\mathcal{A}}, \delta^{\mathcal{B}}, \boldsymbol{X}_\text{N} \right) 
&=\frac{1}{N^2}\sum_{h \in \mathcal{A}}{\sum_{i \in \mathcal{A}}{\pi^{-1}_{h,i}\left\{\frac{\pi_{h,i}}{\pi_{h}\pi_{i}}-1\right\}\left\{\frac{G\left(\min{\left\{R_{h}, R_{i}\right\}}\right)}{n_\text{B}-1}  + G(R_{h})G(R_{i})\right\}}} \\
&+ \frac{1}{N^2}\sum_{h \in \mathcal{A}}{\sum_{i \in \mathcal{A}}{\pi^{-1}_{h,i}\left\{\frac{\pi_{h,i}}{\pi_{h}\pi_{i}}-1\right\}\left\{\frac{-G\left(\min\left\{R_{h}, R_{i}\right\}\right)}{n_\text{B} - 1}\right\}}} \\
&= \frac{1}{N^2}\sum_{h \in \mathcal{A}}{\sum_{i \in \mathcal{A}}{\pi^{-1}_{h,i}\left\{\frac{\pi_{h,i}}{\pi_{h}\pi_{i}}-1\right\}G(R_{h})G(R_{i})}}.
\numberthis \label{final_piece_of_v1}
\end{align*}
Taking the design-based expectation of \eqref{final_piece_of_v1} yields
\begin{align*}
\frac{1}{N^2}\sum_{u \in \mathcal{U}}{\sum_{v \in \mathcal{U}}{\left\{\frac{\pi_{u,v}}{\pi_{u}\pi_{v}}-1\right\}G(R_{u})G(R_{v})}},
\end{align*} and hence 
\begin{align}
\mathbb{E}_{\mathscr{D}}\left[{\mathbb{E}_{\xi}}\left(\widetilde{V}_{1} \ \big| \ \delta^{\mathcal{A}}, \delta^{\mathcal{B}}, \boldsymbol{X}_\text{N} \right)\right] &= \frac{1}{N^2}\sum_{u \in \mathcal{U}}{\sum_{v \in \mathcal{U}}{\left\{\frac{\pi_{u,v}}{\pi_{u}\pi_{v}}-1\right\}G(R_{u})G(R_{v})}}. \label{proof_v1}
\end{align}
The inner expectation of \eqref{unbaised_varest_2} can be expanded in a similar manner:
\begin{align*}
\mathbb{E}_{\xi}\left(\widetilde{V}_{2} \ \big| \ \delta^{\mathcal{A}}, \delta^{\mathcal{B}}, \boldsymbol{X}_\text{N} \right) &= 
\frac{1}{(n_\text{B} - 1)N^2}\sum_{h \in \mathcal{A}}{\sum_{i \in \mathcal{A}}{\pi^{-1}_{h,i}\left\{\frac{\pi_{h,i}}{\pi_{h}\pi_{i}}\right\}}\left\{G\left(\min\left\{R_{h}, R_{i}\right\}\right)\right\}} \\
&- \frac{1}{(n_\text{B} - 1)N^2}\sum_{h \in \mathcal{A}}{\sum_{i \in \mathcal{A}}{\pi^{-1}_{h,i}\left\{\frac{\pi_{h,i}}{\pi_{h}\pi_{i}}\right\}}\left\{\frac{G\left(\min{\left\{R_{h}, R_{i}\right\}}\right)}{ n_\text{B}}  + \frac{n_\text{B} -1 }{n_\text{B}}G(R_{h})G(R_{i})\right\}} \\
&= \frac{1}{(n_\text{B} - 1)N^2}\sum_{h \in \mathcal{A}}{\sum_{i \in \mathcal{A}}{\pi^{-1}_{h,i}\left\{\frac{\pi_{h,i}}{\pi_{h}\pi_{i}}\right\}}\left\{\frac{n_\text{B} - 1}{n_\text{B}}G\left(\min\left\{R_{h}, R_{i}\right\}\right) - \frac{n_\text{B} -1 }{n_\text{B}}G(R_{h})G(R_{i})\right\}} \\
&= \frac{1}{n_\text{B}N^2}\sum_{h \in \mathcal{A}}{\sum_{i \in \mathcal{A}}{\pi^{-1}_{h,i}\left\{\frac{\pi_{h,i}}{\pi_{h}\pi_{i}}\right\}}\big\{G\left(\min\left\{R_{h}, R_{i}\right\}\right) - G(R_{h})G(R_{i})\big\}}.
\numberthis \label{final_piece_of_v2}
\end{align*} Taking the design-based expectation of \eqref{final_piece_of_v2} yields
$$
\frac{1}{n_\text{B}N^2}\sum_{u\in \mathcal{U}}{\sum_{v \in \mathcal{U}}{\left\{\frac{\pi_{u,v}}{\pi_{u}\pi_{v}}\right\}}\big\{G\left(\min\left\{R_{u}, R_{v}\right\}\right) - G(R_{u})G(R_{v})\big\}},
$$ and thus 
\begin{align}
\mathbb{E}_{\mathscr{D}}\left[{\mathbb{E}_{\xi}}\left(\widetilde{V}_{2} \ \big| \ \delta^{\mathcal{A}}, \delta^{\mathcal{B}}, \boldsymbol{X}_\text{N} \right)\right] = \frac{1}{n_\text{B}N^2}\sum_{u\in \mathcal{U}}{\sum_{v \in \mathcal{U}}{\left\{\frac{\pi_{u,v}}{\pi_{u}\pi_{v}}\right\}}\big\{G\left(\min\left\{R_{u}, R_{v}\right\}\right) - G(R_{u})G(R_{v})\big\}}.\label{proof_v2}
\end{align}

Finally, by adding \eqref{proof_v1} and \eqref{proof_v2}, we obtain 
\begin{align*}
\mathbb{E}\left(\widetilde{V}(t; \boldsymbol{\beta}^{*})\right)  &= \frac{1}{N^2}\sum_{u \in \mathcal{U}}{\sum_{v \in \mathcal{U}}{\left\{\frac{\pi_{u,v}}{\pi_{u}\pi_{v}}-1\right\}\big\{G(R_{u})G(R_{v})\big\}}} \\
&+ \frac{1}{n_\text{B}N^2}\sum_{u\in \mathcal{U}}{\sum_{v \in \mathcal{U}}{\left\{\frac{\pi_{u,v}}{\pi_{u}\pi_{v}}\right\}}\big\{G\left(\min\left\{R_{u}, R_{v}\right\}\right) - G(R_{u})G(R_{v})\big\}} \\
&= \mathrm{Var}\left(\widehat{F}_\text{R}(t; \boldsymbol{\beta}^{*})\right),
\end{align*}
which concludes the proof.
\end{proof}

\begin{proof}[Proof of Theorem \ref{theorem_4}]
$\widehat{V}(t; \boldsymbol{\widehat{\beta}}) - \widetilde{V}(t; \boldsymbol{\beta}^{*})$ can be decomposed as the sum of 
\begin{align*}
Q_{1} &= \frac{n_\text{B}}{(n_\text{B} - 1)N^2}\sum_{h \in \mathcal{A}}{\sum_{i \in \mathcal{A}}{\pi^{-1}_{h,i}\left\{\frac{\pi_{h,i}}{\pi_{h}\pi_{i}}-1\right\}\left\{\widehat{G}\big(\widehat{R}_{h}\big)\widehat{G}\big(\widehat{R}_{i}\big) - \widehat{G}\big({R}_{h}\big)\widehat{G}\big({R}_{i}\big)\right\}}}  \\
Q_{2} &=\frac{-1}{(n_\text{B} - 1)N^2}\sum_{h \in \mathcal{A}}{\sum_{i \in \mathcal{A}}{\pi^{-1}_{h,i}\left\{\frac{\pi_{h,i}}{\pi_{h}\pi_{i}}\right\}\left\{\widehat{G}\big(\widehat{R}_{h}\big)\widehat{G}\big(\widehat{R}_{i}\big) - \widehat{G}\big({R}_{h}\big)\widehat{G}\big({R}_{i}\big)\right\}}}  \\
Q_{3} &= \frac{1}{(n_\text{B}-1)N^2}\sum_{h \in \mathcal{A}}{\sum_{i \in \mathcal{A}}{\pi^{-1}_{h,i}}\left\{\widehat{G}\left(\min\left\{\widehat{R}_{h}, \widehat{R}_{i}\right\}\right) - \widehat{G}\left(\min\left\{{R}_{h}, {R}_{i}\right\}\right)\right\}}.\label{vcons_3}
\end{align*} Hence, Theorem \ref{theorem_4} follows by showing that each of $Q_{1}$, $Q_{2}$, and $Q_{3}$ converges in probability to 0.

We first note that $Q_{1}$ can be rewritten as
\begin{align*}
Q_{1} &= \frac{n_\text{B}}{n_\text{B}-1} \times \frac{1}{N^2}\sum_{h \in \mathcal{A}}{\sum_{i \in \mathcal{A}}{\Delta_{1, \pi} \times \pi^{-1}_{h}\pi^{-1}_{i} \left\{\widehat{G}\big(\widehat{R}_{h}\big)\widehat{G}\big(\widehat{R}_{i}\big) - \widehat{G}\big({R}_{h}\big)\widehat{G}\big({R}_{i}\big)\right\}}},\\
&= \frac{1}{1-\frac{1}{n_\text{B}}} \times \frac{1}{N^2}\sum_{h \in \mathcal{A}}{\sum_{i \in \mathcal{A}}{\Delta_{1, \pi} \times \pi^{-1}_{h}\pi^{-1}_{i} \left\{\widehat{G}\big(\widehat{R}_{h}\big)\widehat{G}\big(\widehat{R}_{i}\big) - \widehat{G}\big({R}_{h}\big)\widehat{G}\big({R}_{i}\big)\right\}}},
\end{align*} where $\Delta_{1,\pi} \coloneqq \frac{\pi_{h,i} - \pi_{h}\pi_{i}}{\pi_{h,i}}$. By assumption, $\frac{1}{1-\frac{1}{n_\text{B}}} \leq 2$ and $$-\infty < M_{1} \leq \frac{\pi_{h,i} - \pi_{h}\pi_{i}}{\pi_{h,i}} \leq M_{2} < \infty.$$ 
Because
$$
\widehat{F}^{2}_\text{R}(t; \boldsymbol{\widehat{\beta}}) - \widehat{F}^{2}_\text{R}(t; \boldsymbol{\beta^{*}}) = \frac{1}{N^2}\sum_{h \in \mathcal{A}}{\sum_{i \in \mathcal{A}}}{\pi^{-1}_{h}\pi^{-1}_{i}}\left\{\widehat{G}\big(\widehat{R}_{h}\big)\widehat{G}\big(\widehat{R}_{i}\big) - \widehat{G}\big({R}_{h}\big)\widehat{G}\big({R}_{i}\big)\right\},
$$ it follows that \begin{align*}
&\min\left\{2M_{1} \times \left[\widehat{F}^{2}_\text{R}(t; \boldsymbol{\widehat{\beta}}) - \widehat{F}^{2}_\text{R}(t; \boldsymbol{\beta^{*}}) \right], \ 2M_{2} \times \left[\widehat{F}^{2}_\text{R}(t; \boldsymbol{\widehat{\beta}}) - \widehat{F}^{2}_\text{R}(t; \boldsymbol{\beta^{*}}) \right] \right\} \\
&\leq Q_{1} \leq \\
&\max\left\{2M_{1} \times \left[\widehat{F}^{2}_\text{R}(t; \boldsymbol{\widehat{\beta}}) - \widehat{F}^{2}_\text{R}(t; \boldsymbol{\beta^{*}}) \right], \ 2M_{2} \times \left[\widehat{F}^{2}_\text{R}(t; \boldsymbol{\widehat{\beta}}) - \widehat{F}^{2}_\text{R}(t; \boldsymbol{\beta^{*}}) \right] \right\}.
\end{align*} From the proof of Theorem \ref{theorem_2},  
\begin{align*}
2M_{1} \times \left\{\widehat{F}^{2}_\text{R}(t; \boldsymbol{\widehat{\beta}}) - \widehat{F}^{2}_\text{R}(t; \boldsymbol{\beta^{*}}) \right\} &\xrightarrow{P} 0  \\
2M_{2} \times \left\{\widehat{F}^{2}_\text{R}(t; \boldsymbol{\widehat{\beta}}) - \widehat{F}^{2}_\text{R}(t; \boldsymbol{\beta^{*}}) \right\} &\xrightarrow{P} 0;
\end{align*} and since \begin{align*}
\left\{\left|Q_{1}\right| \geq \varepsilon \right\} &\subseteq \left\{\left|2M_{1} \times \left\{\widehat{F}^{2}_\text{R}(t; \boldsymbol{\widehat{\beta}}) - \widehat{F}^{2}_\text{R}(t; \boldsymbol{\beta^{*}}) \right\}\right| \geq \varepsilon \right\} \\ 
&\cup \left\{\left|2M_{2} \times \left\{\widehat{F}^{2}_\text{R}(t; \boldsymbol{\widehat{\beta}}) - \widehat{F}^{2}_\text{R}(t; \boldsymbol{\beta^{*}}) \right\}\right| \geq \varepsilon \right\} \\
&\coloneqq \left\{\mathcal{A}\right\} \ \cup \ \left\{\mathcal{B}\right\},
\end{align*} we have 
\begin{align*}
\Pr\left\{\left|Q_{1}\right| \geq \varepsilon \right\} &\leq \Pr\left\{\mathcal{A}\right\} + \Pr\left\{\mathcal{B}\right\}  \\
&= 0
\end{align*}  for every $\varepsilon > 0$. Thus $Q_{1} \xrightarrow{P} 0.$ 

Similarly, $Q_{2}$ can be rewritten as
\begin{align*}
Q_{2} &= \frac{-1}{n_\text{B}-1} \times \frac{1}{N^2}\sum_{h \in \mathcal{A}}{\sum_{i \in \mathcal{A}}{\Delta_{2,\pi} \times \pi^{-1}_{h}\pi^{-1}_{i} \left\{\widehat{G}\big(\widehat{R}_{h}\big)\widehat{G}\big(\widehat{R}_{i}\big) - \widehat{G}\big({R}_{h}\big)\widehat{G}\big({R}_{i}\big)\right\}}},
\end{align*}
where 
\begin{align*}
\Delta_{2, \pi} &= \pi^{-1}_{h,i} \times \left\{\frac{\pi_{h,i}}{\pi_{h}\pi_{i}}\right\} \\
&= \pi^{-1}_{h}\pi^{-1}_{i}
\end{align*} is bounded above by some $c \in \mathbb{R}^{+}$ under Assumption \ref{a4A}.  Now, because $-1 \leq \frac{-1}{n_\text{B} -1} \leq  1$, it follows that
\begin{align*}
&\min\left\{-c \times \left[\widehat{F}^{2}_\text{R}(t; \boldsymbol{\widehat{\beta}}) - \widehat{F}^{2}_\text{R}(t; \boldsymbol{\beta^{*}}) \right], c\times \left[\widehat{F}^{2}_\text{R}(t; \boldsymbol{\widehat{\beta}}) - \widehat{F}^{2}_\text{R}(t; \boldsymbol{\beta^{*}}) \right] \right\} \\
&\leq Q_{2} \leq \\
&\max\left\{-c \times \left[\widehat{F}^{2}_\text{R}(t; \boldsymbol{\widehat{\beta}}) - \widehat{F}^{2}_\text{R}(t; \boldsymbol{\beta^{*}}) \right], c \times \left[\widehat{F}^{2}_\text{R}(t; \boldsymbol{\widehat{\beta}}) - \widehat{F}^{2}_\text{R}(t; \boldsymbol{\beta^{*}}) \right] \right\},
\end{align*} and thus $Q_{2} \xrightarrow{P} 0 $.

It now remains to show $Q_{3} \xrightarrow{P} 0$. Since 
$$
\Pr\left(\left|Q_{3}\right| \geq \varepsilon\right) \leq \frac{1}{\varepsilon}\times \mathbb{E}\left|Q_{3}\right|
$$ by Markov's inequality, and 
\begin{align*}
\left|Q_{3}\right|  &\leq \frac{1}{(n_\text{B}-1)N^2}\sum_{h \in \mathcal{A}}{\sum_{i \in \mathcal{A}}}{\pi^{-1}_{h,i}\left|\widehat{G}\left(\min\left\{\widehat{R}_{h}, \widehat{R}_{i}\right\}\right) - \widehat{G}\left(\min\left\{{R}_{h}, {R}_{i}\right\}\right)\right|} \\
&\leq \frac{1}{(n_\text{B}-1)N^2}\sum_{h \in \mathcal{A}}{\sum_{i \in \mathcal{A}}}{\pi^{-1}_{h, i}} \\
&\leq \frac{1}{n_\text{B}-1}  \times M_{2} \times \left(\frac{n_\text{A}}{N}\right)^2,
\end{align*} we see that $\Pr\left(\left|Q_{3}\right| \geq \varepsilon\right) \to 0$ by the monotonicity property of expectations. Thus $Q_{3} \xrightarrow{P} 0$, which concludes the proof.
\end{proof}

\end{document}